

%
%
%
%
\ifx\UsualIsLoaded\undefined
\let\UsualIsLoaded=\relax		
\font\fourteenrm=cmr12  scaled \magstep1
\font\fourteenbf=cmbx12 scaled \magstep1
\font\fourteentt=cmtt12 scaled \magstep1
\font\fourteensl=cmsl12 scaled \magstep1
\font\fourteensy=cmsy10 scaled \magstep2
\font\fourteeni=cmmi12  scaled \magstep1
\font\fourteenit=cmti12 scaled \magstep1
\font\fourteensc=cmcsc10 at 14pt 
\font\fourteenbit=cmssi12 scaled \magstep1
\font\fourteenbbb=msbm10 scaled \magstep2
\font\twelverm=cmr12
\font\twelvebf=cmbx12
\font\twelvett=cmtt12
\font\twelvesl=cmsl12
\font\twelvesy=cmsy10 scaled \magstep1
\font\twelvei=cmmi12
\font\twelveit=cmti12
\font\twelvesc=cmcsc10 at 11.4pt 
\font\twelvebit=cmssi12
\font\twelvebbb=msbm10 scaled \magstep1
\font\tenrm=cmr10
\font\tenbf=cmb10 
\font\tentt=cmtt10 
\font\tensl=cmsl10 
\font\tensy=cmsy10 
\font\teni=cmmi10 
\font\tenit=cmti10 
\font\tensc=cmcsc10 at 9.5pt
\font\tenbit=cmssi10
\font\tenbbb=msbm10
\font\ninei=cmmi9 
\font\ninerm=cmr9

\font\ninesy=cmsy9 
\font\ninebbb=msbm9
\font\eighti=cmmi8
\font\eightrm=cmr8

\font\eightsy=cmsy8
\font\eightbbb=msbm8
\font\seveni=cmmi7 
\font\sevenrm=cmr7
\font\sevensy=cmsy7
\font\sevenbbb=msbm7
%
%
%

%

%
\newfam\bbbfam
%
\def\tenpoint{
\def\rm{\fam0\tenrm}%
\textfont0=\tenrm \scriptfont0=\eightrm \scriptscriptfont0=\sevenrm
\textfont1=\teni \scriptfont1=\eighti \scriptscriptfont1=\seveni
\textfont2=\tensy \scriptfont2=\eightsy \scriptscriptfont2=\sevensy
\textfont3=\tenex \scriptfont3=\tenex \scriptscriptfont3=\tenex
\textfont\itfam=\tenit \def\it{\fam\itfam\tenit}%
\textfont\slfam=\tensl \def\sl{\fam\slfam\tensl}%
\textfont\bffam=\tenbf \def\bf{\fam\bffam\tenbf}%
\textfont\ttfam=\tentt \def\tt{\fam\ttfam\tentt}%
\def\sc{\tensc}%
\def\bit{\tenbit}%
\textfont\bbbfam=\tenbbb \scriptfont\bbbfam=\eightbbb
\scriptscriptfont\bbbfam=\sevenbbb \def\bbb{\fam\bbbfam\twelvebbb}%
\normalbaselineskip=12pt
\setbox\strutbox=\hbox{\vrule height10pt depth4pt width0pt}%
\normalbaselines\rm}
%
%
\def\twelvepoint{
\def\rm{\fam0\twelverm}%
\textfont0=\twelverm \scriptfont0=\ninerm \scriptscriptfont0=\sevenrm
\textfont1=\twelvei \scriptfont1=\ninei \scriptscriptfont1=\seveni
\textfont2=\twelvesy \scriptfont2=\ninesy \scriptscriptfont2=\sevensy
\textfont3=\tenex \scriptfont3=\tenex \scriptscriptfont3=\tenex
\textfont\itfam=\twelveit  \def\it{\fam\itfam\twelveit}%
\textfont\slfam=\twelvesl  \def\sl{\fam\slfam\twelvesl}%
\textfont\bffam=\twelvebf  \def\bf{\fam\bffam\twelvebf}%
\textfont\ttfam=\twelvett  \def\tt{\fam\ttfam\twelvett}%
\def\sc{\twelvesc}%
\def\bit{\twelvebit}%
\def\bbb{\fam\bbbfam\twelvebbb}%
\textfont\bbbfam=\twelvebbb \scriptfont\bbbfam=\ninebbb
\scriptscriptfont\bbbfam=\sevenbbb
\normalbaselineskip=14pt
\setbox\strutbox=\hbox{\vrule height10pt depth4pt width0pt}%
\normalbaselines\rm}
%
%
\def\fourteenpoint{
\def\rm{\fam0\fourteenrm}%
\textfont0=\fourteenrm \scriptfont0=\twelverm \scriptscriptfont0=\tenrm
\textfont1=\fourteeni \scriptfont1=\twelvei \scriptscriptfont1=\teni
\textfont2=\fourteensy \scriptfont2=\twelvesy \scriptscriptfont2=\tensy
\textfont3=\tenex \scriptfont3=\tenex \scriptscriptfont3=\tenex
\textfont\itfam=\fourteenit \def\it{\fam\itfam\fourteenit}%
\textfont\slfam=\fourteensl \def\sl{\fam\slfam\fourteensl}%
\textfont\bffam=\fourteenbf \def\bf{\fam\bffam\fourteenbf}%
\textfont\ttfam=\fourteentt \def\tt{\fam\ttfam\fourteentt}%
\def\sc{\fourteensc}%
\def\bit{\fourteenbit}%
\def\bbb{\fam\bbbfam\twelvebbb}%
\textfont\bbbfam=\fourteenbbb \scriptfont\bbbfam=\tenrm 
\scriptscriptfont\bbbfam=\eightrm
\normalbaselineskip=16pt
\setbox\strutbox=\hbox{\vrule height10pt depth4pt width0pt}%
\normalbaselines\rm}
%
%
%
\twelvepoint
\abovedisplayskip 14pt plus 3pt minus 10pt%
\belowdisplayskip 14pt plus 3pt minus 10pt%
\abovedisplayshortskip 0pt plus 3pt%
\belowdisplayshortskip 8pt plus 3pt minus 5pt%
\parskip 3pt plus 1.5pt
\hsize=6.5in
\vsize=8.9in
\fi			


\input psfig.sty
\input amssym.def
\input amssym
\def\widehat{\mathaccent"0362 }
\def\ss{\smallskip}
\def\ms{\medskip}
\def\bs{\bigskip}
\def\ni{\noindent}
\def\cl{\centerline}

\def\QED {$\quad\square$}
\def\QP{\smallskip\leftskip=.4in\rightskip=.4in\noindent}
\def\ref{\hangindent=1pc \hangafter=1 \noindent}

\def\ssm{\smallsetminus}

\def\mod{{\rm mod\;}}

\def\[{$\,}
\def\]{\,$}
\def\th{\kern 3pt}
\def\C{{\bbb C}}

\def\Z{{\bbb Z}}

\def\Q{{\bbb Q}}

\def\>{\succ}

\def\={~=~ }

\def\FIG#1#2{\cl{\psfig{figure=#1,height=#2}}}
\def\tp{\tenpoint\hskip -.1in}

\font\bit=cmssi12 at 12truept
\font\title=cmbx12

\mathsurround = 1pt
\abovedisplayskip=6pt
\belowdisplayskip=6pt
\parskip=2pt

\def\IMSmarkvadjust{0 pt}
\def\IMSmarkhadjust{0 pt}
\def\IMSmarkhpadding{0 pt}
\def\IMSpubltext{Published in modified form:}
\def\SBIMSMark#1#2#3{
 \font\SBF=cmss10 at 10 true pt
 \font\SBI=cmssi10 at 10 true pt
 \setbox0=\hbox{\SBF \hbox to \IMSmarkhpadding{\relax}
                Stony Brook IMS Preprint \##1}
 \setbox2=\hbox to \wd0{\hfil \SBI #2}
 \setbox4=\hbox to \wd0{\hfil \SBI #3}
 \setbox6=\hbox to \wd0{\hss
             \vbox{\hsize=\wd0 \parskip=0pt \baselineskip=10 true pt
                   \copy0 \break%
                   \copy2 \break%
                   \copy4 \break}}
 \dimen0=\ht6   \advance\dimen0 by \vsize \advance\dimen0 by 8 true pt
                \advance\dimen0 by -\pagetotal
	        \advance\dimen0 by \IMSmarkvadjust
 \dimen2=\hsize \advance\dimen2 by .25 true in
	        \advance\dimen2 by \IMSmarkhadjust

%
%
  \openin2=publishd.tex
  \ifeof2\setbox0=\hbox to 0pt{}
  \else 
     \setbox0=\hbox to 3.1 true in{
                \vbox to \ht6{\hsize=3 true in \parskip=0pt  \noindent  
                {\SBI \IMSpubltext}\hfil\break
                \input publishd.tex 
                \vfill}}
  \fi
  \closein2
  \ht0=0pt \dp0=0pt
 \ht6=0pt \dp6=0pt
 \setbox8=\vbox to \dimen0{\vfill \hbox to \dimen2{\copy0 \hss \copy6}}
 \ht8=0pt \dp8=0pt \wd8=0pt
 \copy8
 \message{*** Stony Brook IMS Preprint #1, #2. #3 ***}
}

\def\T{{\cal T}}
\def\H{{\bbb H}}
\def\P{{\bbb P}}
\def\d{{\bf d}}
\def\cn{{\rm cn}}
\def\dn{{\rm dn}}
\def\sn{{\rm sn}}
\def\wC{{\widehat\C}}
\def\mapsdown{\psfig{figure=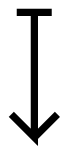,height=.11in}}
\def\mapsfrom{\psfig{figure=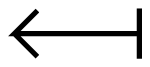,height=.09in}}

\def\lloop{\psfig{figure=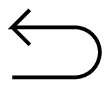,height=.11in}}
\def\mapsup{\psfig{figure=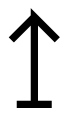,height=.11in}}
\def\obul{{\bullet\hskip -.106in\odot}}
\font\tencyr=wncyr10
\input cyracc.def
\def\cyr{\tencyr\cyracc}
\def\Ch{\hbox{\cyr Q}}

\SBIMSMark{2004/01}{Feb.~2004, revised Sept.~2004}{}

\footline={}
\def\newsech#1#2{\count255=\pageno
        \headline={%
    \ifnum\pageno > \count255 {%
        \hss\tenrm\ifodd\pageno #1\hss\folio
        \else\folio\hfill #2\hfill\fi}
    \else\ifodd\pageno\hss\tenrm\folio
        \else\tenrm\folio\hfill
        \fi
    \fi}}
\newsech{\sc On Latt\`es Maps}{\sc John Milnor}

\cl{\bf On Latt\`es Maps}\ss

\cl{John Milnor}\ss

\cl{Dedicated to Bodil Branner.}\ms

{\leftskip=.8in\rightskip=.8in\ni
\bit Abstract. An exposition of the 1918 paper of Latt\`es, together with
its historical antecedents, and its modern formulations and applications.\ss}

\hskip 1.6in 1. The Latt\`es paper.

\hskip 1.6in 2. Finite Quotients of Affine Maps

\hskip 1.6in 3. A Cyclic Group Action on \[\C/\Lambda\].

\hskip 1.6in 4. Flat Orbifold Metrics

\hskip 1.6in 5. Classification

\hskip 1.6in 6. Latt\`es Maps before Latt\`es

\hskip 1.6in 7. More Recent Developments

\hskip 1.6in 8. Examples

\hskip 1.6in References\ms

{\bf\S1. The Latt\`es paper.}
In 1918, some months before his death of typhoid fever, Samuel Latt\`es
published a brief paper describing an extremely interesting class of rational
maps. Similar examples had been described by Schr\"oder almost fifty
years earlier (see \S6), but Latt\`es'
name has become firmly attached to these maps, which play a basic role
as exceptional examples in the holomorphic dynamics literature.

His starting point
was the ``{\bit Poincar\'e function\/}'' \[\theta:\C~\to~\widehat\C\]
associated with a repelling fixed point \[z_0=f(z_0)\] of a rational function
\[f:\widehat\C\to\widehat\C\]. This can be described as the
inverse of  the K{\oe}nigs linearization around \[z_0\],
extended to a globally defined meromorphic function.\footnote{$^1$}{\tp
Compare [La], [P], [K].
For general background material, see for example [M3] or [BM].} Assuming for
convenience that \[z_0\ne\infty\], it is characterized by the identity
$$  f(\theta(t))\=\theta(\mu\,t) $$
for all complex numbers \[t\], with \[\theta(0)=z_0\],  normalized by
the condition that \[\theta'(0)=1\]. Here \[\mu=f'(z_0)\] is the
{\bit multiplier\/} at \[z_0\], with \[|\mu|>1\]. This Poincar\'e function
can be computed explicitly by the formula
$$	\theta(t)\=\lim_{n\to\infty}\,f^{\circ n}\Big(z_0+t/\mu^n\Big)~. $$
Its image \[\theta(\C)\subset\wC\] is equal to the Riemann sphere \[\wC\]
with at most two points removed. In practice, we will always assume that \[f\]
has degree at least two. The complement \[\wC\ssm\theta(\C)\] is then
precisely equal to the {\bit exceptional set\/} \[{\cal E}_f\],
consisting of all points with finite grand orbit under \[f\].

In general this Poincar\'e
function \[\theta\] has very complicated behavior.
In particular, the Poincar\'e functions associated with different fixed
points or periodic points are usually
quite incompatible. However, Latt\`es pointed
out that in special cases \[\theta\] will be periodic or doubly periodic,
and will give rise to a simultaneous linearization
for all of the periodic points of \[f\].
(For a more precise statement, see the proof of 3.9 below.)

We will expand on this idea in the following sections.
Section 2 will introduce rational maps which are
{\bit finite quotients of affine maps\/}. (These are more
commonly described in the literature as rational maps with {\bit flat orbifold
metric\/}---see \S4.) They can be classified into {\bit power maps\/},
{\bit Chebyshev maps\/}, and {\bit Latt\`es maps\/}
according as the Julia set is a circle,
a line or circle segment, or the entire Riemann sphere. These
maps will be studied in Sections 3 through 5, concentrating on the Latt\`es
case. Section 6
will describe the history of these ideas before Latt\`es; and \S7
will describe some of the developments since his time.
Finally, \S8 will describe a number of concrete examples.\ss

{\bf Acknowledgments.} I want to thank Curt McMullen, Alexandre
Eremenko, and Walter Bergweiler for their help, both with the mathematics
and with the history, and I want to thank the National Science Foundation
(DMS 0103646) and the
Clay Mathematics Institute for their support of Dynamical System activities
in Stony Brook.\bs

{\bf \S2. Finite Quotients of Affine Maps.}
It will be convenient to make a very mild generalization of the Latt\`es
construction, replacing the linear map \[t\mapsto \mu\,t\] of his construction
by an affine map \[t\mapsto a\,t+b\]. Let \[\Lambda\] be a discrete additive
subgroup of the complex numbers \[\C\].
In the cases of interest, this subgroup will have rank either
one or two, so that the quotient surface \[\C/\Lambda\] is either a cylinder
\[\cal C\] or a torus \[\cal T\].
\ss

{\bf Definition 2.1.} A rational map \[f\] of degree two or more will be
called a {\bit finite quotient of an affine map\/} if there is a flat surface
\[\C/\Lambda\], an affine map \[L(t)=a\,t+b\] from \[\C/\Lambda\] to itself,
and  a finite-to-one holomorphic map
\[\Theta:\C/\Lambda\to\wC\ssm {\cal E}_f\] which satisfies
the semiconjugacy relation \[f\circ\Theta=\Theta\circ L\].
Thus the following diagram must commute:
$$\matrix{~\C/\Lambda & \buildrel L\over\longrightarrow & ~\C/\Lambda\cr
\Theta\downarrow~ &&\Theta\downarrow~\cr
~~\wC\ssm {\cal E}_f & \buildrel f\over\longrightarrow &
~~~\wC\ssm {\cal E}_f\,.} \eqno(1)$$
We can also write \[f=\Theta\circ L\circ\Theta^{-1}\].
It follows for example that any periodic orbit of \[L\] must map to a
periodic orbit of \[f\], and conversely that every periodic orbit of \[f\]
outside of the exceptional set \[{\cal E}_f\] is
the image of a periodic orbit of \[L\]. (However, the periods are not
necessarily the same.) Here the finite-to-one condition is essential. In
fact it follows from Poincar\'e's construction that any rational map of
degree at least two can be thought of as an
infinite-to-one quotient of an affine map of \[\C\].

These finite quotients of affine maps can be classified very roughly into three
types, as follows.
The set of postcritical points of \[f\] plays an important role in all cases.
(Compare Lemma 3.4.)\ss

{\bf Power maps.} These are the simplest examples.
By definition, a rational map will be called a
{\bit power map\/} if it is holomorphically conjugate to a map of the form
$$    f_a(z)~\=~ z^a  $$
where \[a\] is an integer. Note that \[f_a\], restricted to the punctured plane
\[\C\ssm\{0\}=\wC\ssm{\cal E}_{f_a}\], is conjugate to the linear map
\[t\mapsto at\] on the cylinder \[\C/2\pi\Z\]. In fact
\[f_a(e^{it})=e^{i\,at}\],
where the conjugacy \[t\mapsto e^{it}\] maps
\[\C/2\pi\Z\] diffeomorphically onto \[\C\ssm\{0\}\]. The degree of \[f_a\] is
equal to \[|a|\], the Julia set \[J(f_a)\] is equal to the unit circle, and the
exceptional set \[{\cal E}_{f_a}=\{0\,,\,\infty\}\] consists of the two
critical points, which are also the two postcritical points.\ss

{\bf Chebyshev maps.} These are the next simplest examples. A rational
map will be called a {\bit Chebyshev map\/} if it is
conjugate to \[\pm\Ch_n(z)\] where \[\Ch_n\] is the
degree \[n\] {\bit Chebyshev polynomial\/},
defined by the equation\footnote{$^2$}{\tp The Russian letter
\[\Ch\] is called ``chi'', pronounced as in ``chicken''.}
$$ \Ch_n(u+u^{-1})\=u^n+u^{-n}~. $$
For example:
$$	\Ch_2(z)\=z^2\,-\,2~,\qquad \Ch_3(z)\= z^3\,-\,3z~,\quad
	\Ch_4(z)\= z^4\,-\,4z^2\,+\,2~,\quad\cdots~.$$
We will see in \S3.8 that power maps and Chebyshev maps
are the only finite quotients of affine maps
for which the lattice \[\Lambda\subset\C\] has rank one.

If we set \[u=e^{it}\], then the map \[\Theta(t)=u+u^{-1}=2\,\cos(t)\]
is a proper map of degree two from the cylinder \[\C/2\pi\Z\] to the plane
\[\C\], satisfying
$$  \Ch_n(\Theta(t))\=\Theta(n\,t)\qquad{\rm or~equivalently}\qquad
\Ch_n(2\,\cos t)\=2\,\cos(nt)~,$$
and also \[-\Ch_n(\Theta(t))=\Theta(n\,t+\pi)\]. These identities
show that both \[\Ch_n\] and \[-\Ch_n\]
are finite quotients of affine maps.
The Julia set \[J(\pm\Ch_n)\] is the closed interval \[[-2\,,\,2]\], and the
exceptional set for \[\pm\Ch_n\] is the singleton \[\{\infty\}\].
The postcritical set of \[\pm\Ch_n\] consists of the three points
\[\{\pm 2\,,\,\infty\}\]. In fact,
if \[2 \cos(t)\] is a finite critical point of \[\Ch_n\] then
by differentiating the equation \[\Ch_n(2\cos\,t)=2 \cos(n\,t)\] we see
that \[\sin(n\,t)=0\] and hence that \[2\cos(n\,t)=\pm 2\].

Note: If \[n\] is even, the equation
\[-\Ch_n(z)=\Ch_n(kz)/k\] with \[k=-1\] shows that \[-\Ch_n\] is holomorphically conjugate to
\[\Ch_n\]. However, for \[n\] odd the map \[z\mapsto -\Ch_n(z)\] has a
postcritical orbit \[\{\pm 2\}\] of period two, and hence cannot be conjugate
to \[z\mapsto\Ch_n(z)\] which has only postcritical fixed points.\ss

{\bf Latt\`es maps.} In the remaining case where the lattice \[\Lambda
\subset\C\] has rank two so that the quotient \[\T=\C/\Lambda\] is a torus,
the map \[f=\Theta\circ L\circ\Theta^{-1}\] will be called
a {\bit Latt\`es map\/}.
Here \[L\] is to be an affine self-map of the torus,
and \[\Theta\] is to be a holomorphic map from \[\T\]
to the Riemann sphere \[\widehat\C\]. These are the most
interesting examples, and exhibit rather varied behavior. Thus we can
distinguish between {\bit flexible\/} Latt\`es maps
which admit smooth deformations, and rigid Latt\`es maps which do not.
(See 5.5 and 5.6, as well as \S7 and 8.3.)
Another important distinction is between the  Latt\`es maps
with three postcritical points, associated with triangle groups acting on
the plane, and those with four postcritical points. (See \S4.)

For any Latt\`es map \[f\], since \[\Theta\] is necessarily onto, there are
no exceptional points. Furthermore, since periodic points
of \[L\] are dense on the torus it follows that periodic points of \[f\]
are dense on the Riemann sphere. Thus the Julia set \[J(f)\]
must be the entire sphere.\ms

{\bf\S3 Cyclic Group Actions on \[\C/\Lambda\].} The following
result provides a more explicit description of all of the possible
Latt\`es maps, as defined in \S2.

{\QP{\bf Theorem 3.1.} \it A rational map is Latt\`es
if and only if it is conformally conjugate to a map of the form
$\;	L/G_n\,:\,\T/G_n~\to~\T/G_n \;$  where:

\ni\[\bullet\] \[\T\cong\C/\Lambda\] is a flat torus,

\ni\[\bullet\] \[G_n\] is the group of \[n$-th roots of unity acting on
\[\T\] by rotation around a base point, with \[n\] equal to
either~ 2, 3, 4, or 6,

\ni\[\bullet\]
\[\T/G_n\] is the quotient space provided with its natural structure
as a smooth Riemann surface of genus zero,

\ni\[\bullet\] \[L\] is an affine map from \[\T\] to itself which commutes
with a generator of \[G_n\], and

\ni\[\bullet\] \[L/G_n\] is the induced holomorphic map from
the quotient surface to itself.\ms}

{\bf Remark 3.2.} The map \[\T\to\T/G_n\cong\wC\] can of course
be described in terms of classical elliptic function theory. In the case
\[n=2\] we can identify this map with the Weierstrass function \[\wp:\C/\Lambda
\to\wC\] associated with
the period lattice \[\Lambda\]. Here the lattice \[\Lambda\] or the torus
\[\T\] can be completely
arbitrary, but in the cases \[n\ge 3\] we will see that \[\T\] is uniquely
determined by \[n\], up to conformal isomorphism.
For \[n=3\] we can take the derivative
\[\wp'\] of the associated Weierstrass function as the semiconjugacy
\[\wp':\T\to\wC\], while for \[n=6\] we can use either \[(\wp')^2\]
or \[\wp^3\] as semiconjugacy. (For any lattice with \[G_3$-symmetry, these
two functions are related by the identity \[(\wp')^2=4\wp^3+{\rm constant}\].
The two alternate forms correspond to the fact
that \[\T/G_6\] can be identified either with \[(\T/G_3)/G_2\] or with
\[(\T/G_2)/G_3\].) Finally, for \[n=4\] we can use the square \[\wp^2\]
of the associated Weierstrass function, corresponding to the factorization
\[\T\to\T/G_2\to\T/G_4\].

{\bf Remark 3.3.} This theorem is related to the definition in \S2 as follows.
Let us use the notation \[\Theta^\star:\T^\star=\C/\Lambda^\star\to\wC\ssm{\cal E}_f\]
for the initial semiconjugacy of Definition 2.1, formula (1). The
degree of this semiconjugacy \[\Theta^\star\]
can be arbitrarily large. However, the proof of 3.1 will show that
\[\Theta^\star\] can be factored in an essentially unique way
as a composition \[\T^\star\to\T\to\T/G_n\cong\wC\] for some torus \[\T\],
with \[n\] equal to 2, 3, 4, or 6.\ms

The proof of 3.1 will be based on the following ideas.
Let \[\theta:\C\to\wC\] be a doubly periodic meromorphic function, and let
\[\Lambda\subset\C\] be its lattice of periods so that
\[\lambda\in\Lambda\] if and only if \[\theta(t+\lambda)=\theta(t)\] for all
\[t\in\C\]. Then the canonical flat metric \[|\d t|^2\] on \[\C\] pushes
forward to a corresponding flat metric on the torus \[\T=\C/\Lambda\].
If \[\ell(t)=at+b\] is an affine map of \[\C\] satisfying the identity
\[f\circ\theta=\theta\circ\ell\], then for \[\lambda\in\Lambda\] and \[t\in\C\]
we have
$$ \theta(at+b)\=f(\theta(t))\=f(\theta(t+\lambda))\=\theta(a(t+\lambda)+b)~.$$
It follows that \[a\Lambda\subset\Lambda\].
Therefore the maps \[\ell\] and \[\theta\] on \[\C\] induce
corresponding maps \[L\] and \[\Theta\] on \[\T\], so that we have a
commutative diagram of holomorphic maps
$$\matrix{~\T & \buildrel L\over\longrightarrow & ~\T\cr
\Theta\downarrow~ &&\Theta\downarrow~\cr
~\widehat\C & \buildrel f\over\longrightarrow &~\;\widehat\C\,.} \eqno(2)$$
We will think of \[\T\] as a branched covering of the Riemann sphere with
projection map \[\Theta\]. Since \[L\] carries a small region of
area \[A\] to a region of area \[|a|^2A\], it follows that the map \[L\] has
degree \[|a|^2\]. Using Diagram (2), we see that the degree \[d_f\ge 2\]
of the map \[f\] must also be equal to \[|a|^2\].

One easily derived property is the following. (For a more precise statement,
see 4.5.)
Let \[C_f\] be the set of critical points of \[f\] and let \[V_f=f(C_f)\]
be the set of critical values. Similarly,
let \[V_\Theta=\Theta(C_\Theta)\] be the set
of critical values for the projection map \[\Theta\].

{\QP{\bf Lemma 3.4.} \it Every Latt\`es map \[f\] is postcritically
finite. In fact the postcritical set
$$  P_f\=V_f\cup f(V_f)\cup f^{\circ 2}(V_f)\cup\cdots  $$
for \[f\] is precisely equal to the finite set \[V_\Theta\] consisting of all
critical values for the projection \[\Theta:\T\to \widehat\C\].\ms}

{\bf Proof.} Let \[d_f(z)\] be the local degree of the map \[f\] at a point
\[z\]. Thus
$$	1~\le ~d_f(z)~\le~ d_f~, $$
where \[d_f(z)>1\] if and only if \[z\] is a critical point of \[f\].
Given points \[\tau_j\in \T\] and \[z_j\in\widehat\C\] with\vskip -.3in
$$\matrix{~\tau_1 & \buildrel L\over \mapsto & ~\tau_0\cr
	\Theta\,\mapsdown~ && \Theta\,\mapsdown~\cr
	~z_1&\buildrel f\over\mapsto & ~~z_0~,} $$
since \[L\] has local degree \[d_L(\tau)=1\] everywhere,
it follows that
$$  d_\Theta(\tau_0)\=d_\Theta(\tau_1)\cdot d_f(z_1)~.\eqno(3) $$
Since the maps \[L\] and \[\Theta\] are surjective, it follows that
\[z_0\] is a critical value of \[\Theta\] if and and if it is either a critical
value of \[f\] or has a preimage \[z_1\in f^{-1}(z_0)\] which is a critical
value of \[\Theta\] or both. Thus \[\; V_\Theta\=V_f\cup f(V_\Theta)\],
which implies inductively that
$$ f^{\circ n}(V_f)\subset V_{\Theta}~,
\qquad{\rm hence} \qquad P_f~\subset~ V_\Theta~. $$
On the other hand, if some critical point \[\tau_0\] of \[\Theta\] had image
\[\Theta(\tau_0)\] outside of the postcritical set \[P_f\], then all of the
infinitely many iterated preimages \[\cdots\mapsto \tau_2\mapsto \tau_1\mapsto
\tau_0\] would have the same property. This is impossible, since \[\Theta\]
can have only finitely many critical points.\QED
\ss

We will prove the following preliminary version
of 3.1, with notations as in Diagram (2).

{\QP{\bf Lemma 3.5.} \it If \[f\] is a Latt\`es map, then there is a finite
cyclic group \[G\] of rigid rotations of the
torus \[\T\] about some base point,
so that \[\Theta(\tau')=\Theta(\tau)\] if and only if
\[\tau'=g\,\tau\] for some \[g\in G\].
Thus \[\Theta\] induces a canonical homeomorphism from the quotient
space \[\T/G\] onto the Riemann sphere.\ss}

{\bf Remark 3.6.} Such a quotient \[\T/G\] can be
given two different structures which are distinct, but closely related.
Suppose that a point \[\tau_0\in\T\] is mapped to itself by a
subgroup of \[G\], necessarily cyclic, of order \[r>1\].
Any \[\tau\] close to \[\tau_0\] can be written as
the sum of \[\tau_0\] with a small complex number \[\tau-\tau_0\].
The power \[(\tau-\tau_0)^r\] then serves as a local
uniformizing parameter for \[\T/G\] near \[\tau_0\]. In this way, the
quotient becomes a smooth Riemann surface. On the other hand,
if we want to carry the flat Euclidean structure of \[\T\] over
to \[\T/G\], then the image of \[\tau_0\] must be considered as a singular
``cone point'', as described in the next section. The integer
\[r\], equal to the local degree \[d_\Theta(\tau_0)\], depends only
on the image point \[\Theta(\tau_0)\], and
is called the {\bit ramification index\/} of \[\Theta(\tau_0)\].
\ss

{\bf Proof of 3.5.} Let \[U\] be any simply connected open subset of
\[\wC\ssm P_f=\wC\ssm V_\Theta\]. Then the preimage \[\Theta^{-1}(U)\]
is the union \[U_1\cup\cdots\cup U_n\] of \[n\] disjoint open sets,
each of which projects diffeomorphically onto \[U\],
where \[n=d_\Theta\] is equal to the degree of \[\Theta\].
Let \[\Theta_j:U_j\buildrel\cong\over\longrightarrow U\] be the restriction
of \[\Theta\] to \[U_j\].
We will first prove that each composition
$$	\Theta_k^{-1}\circ\Theta_j~:~U_j~\to~U_k ~,\eqno(4)$$
is an isometry from \[U_j\] onto \[U_k\], using the standard flat metric on
the torus.

Since periodic points of \[f\] are everywhere dense, we can choose
a periodic point \[z_0\in U\]. Now replacing \[f\] by some carefully chosen
iterate, and replacing \[L\] by the corresponding iterate, we may assume
without changing \[\Theta\] that:

{\QP\[\bullet\]  \[z_0\] is actually a fixed point of \[f\], and that\par}

{\QP\[\bullet\] every point in the finite \[L$-invariant
set \[\Theta^{-1}(z_0)\] is either fixed by \[L\], or\break
 \phantom{---} maps directly to a fixed point.\ss}

\ni In other words, each point \[\tau_j=\Theta_j^{-1}(z_0)\] is either
a fixed point of \[L\] or maps to a fixed point. For \[\tau\] close to
\[\tau_j\], evidently the difference \[\tau-\tau_j\] can be identified with
a unique complex number close to zero. Setting
\[L(\tau_j)=\tau_{j'}\], note that the affine map
$$   \tau-\tau_j~~\mapsto~~L(\tau)-\tau_{j'}~~\in~~\C $$
is actually linear, so that \[L(\tau)-\tau_{j'}=\mu\,(\tau-\tau_j)\]
where \[\mu=L'\] is constant.
Similarly, for \[z\] close to \[z_0\], the difference
$$\kappa_j(z)\=\Theta_j^{-1}(z)- \tau_j $$
is well defined as a complex number. For each index \[j\] we will show
that the map\break \[~z~\mapsto~\kappa_j(z)~\in~\C~\]
is a K\oe nigs linearizing map for \[f\] in a neighborhood of \[z_0\]. That is,
$$\kappa_j(f(z))\=\mu\,\kappa_j(z)~, \qquad{\rm with} \qquad
\kappa_j(z_0)\=0~,\eqno(5)$$
where the constant \[\mu=L'\] is necessarily equal to
the multiplier of \[f\] at \[z_0\].
In fact the identity \[\Theta_{j'}^{-1}(f(z))=L(\Theta_j^{-1}(z))\] holds
for all \[z\] close to \[z_0\]. Subtracting \[\tau_{j'}\] we see that
$$\kappa_{j'}(f(z))\=\mu\, \kappa_j(z)~.\eqno(6) $$
If \[\tau_j\] is a fixed point
so that \[j'=j\], then this is the required assertion (5). But \[\tau_{j'}\]
is always a fixed point, so this proves that
\[\kappa_{j'}(f(z))\=\mu\,\kappa_{j'}(z)\].  Combining
this equation with (6), we see that \[\kappa_j(z)=\kappa_{j'}(z)\], and it
follows that equation (5) holds in all cases.

Since such a K\oe nigs linearizing map is unique up to multiplication
by a constant, it follows that every \[\kappa_i(z)\] must be equal to the
product \[c_{ij}\,\kappa_j(z)\] for some
constant \[c_{ij}\ne 0\] and for all \[z\] close to \[z_0\]. Therefore
\[\Theta_i^{-1}(z)\] must be equal to \[c_{ij}\,\Theta_j^{-1}(z)\] plus a
constant for all \[z\] close to \[z_0\].
 Choosing a local lifting of
\[\Theta_i^{-1}\circ\Theta_j\] to the universal covering space
\[\widetilde\T\cong\C\] and continuing analytically, we obtain an affine map
\[A_{ij}\] from \[\C\] to itself with derivative \[A_{ij}'=c_{ij}\],
satisfying the identity \[\theta=\theta\circ A_{ij}\], where \[\theta\]
is the composition \[\widetilde\T\to\T\buildrel\Theta\over\longrightarrow\wC\].

We must prove that \[|c_{jk}|=1\], so that this affine transformation
is an isometry. Let \[\widetilde G\] be the group\footnote{$^3$}{\tp
This \[\widetilde G\] is often described as a {\bit crystallographic group\/}
acting on \[\C\].
That is, it is a discrete group of rigid Euclidean motions of \[\C\],
with compact quotient \[\C/\widetilde G\cong\T/G\].}
consisting of all affine transformations \[\widetilde g\] of \[\C\]
which satisfy the identity \[\theta=\theta\circ\widetilde g\].
 The
translations \[t\mapsto t+\lambda\] with \[\lambda\in\Lambda\] constitute
a normal subgroup, and the quotient \[G=\widetilde G/\Lambda\]
acts as a finite group of complex affine automorphisms
of the torus \[\T=\C/\Lambda\].
In fact \[G\] has exactly \[n\] elements, since it contains exactly
one transformation \[g\] carrying \[U_1\] to any specified \[U_j\].
The derivative map
\[g\mapsto g'\] is an injective homomorphism from \[G\] to the
multiplicative group \[\C\ssm\{0\}\]. Hence it must carry \[G\] isomorphically
onto the unique subgroup of \[\C\ssm\{0\}\] of order \[n\], namely the
group\[G_n\] of \[n$-th roots of unity.
Furthermore, a generator of \[G\] must have a
fixed point in the torus, so \[G\] can be considered as a group of rotations
about this fixed point. This completes the proof of 3.5.\QED
\ss

In fact, if we translate coordinates so that some specified fixed point
of the \[G$-action is the origin of the torus
\[\T=\C/\Lambda\], then clearly we can identify \[G\] with the group \[G_n\]
of \[n$-th roots of unity, acting by multiplication on \[\T\].\ss

{\QP{\bf Lemma 3.7.} \it The order \[n\] of such a cyclic group of rotations
of the torus with quotient \[\T/G_n\cong\wC\]
is necessarily either \[2,~3,~4\], or \[6\].\ss}

{\bf Proof.} Thinking of a rotation through angle \[\alpha\] as a real
linear map, it has eigenvalues \[e^{\pm i\,\alpha}\] and trace
\[e^{i\,\alpha}+e^{-i\,\alpha}=2 \cos(\alpha)\]. On the other hand, if such a
rotation
carries the lattice \[\Lambda\] into itself, then its trace must be an integer.
The function \[\alpha\mapsto 2 \cos(\alpha)\] is monotone decreasing for
\[0<\alpha\le\pi\] and takes only the following integer values:
$$\matrix{ r~~&= && 6 & 4 & 3 & 2\cr
2\cos(2\pi/r)&= && 1 & 0 & -1 & -2~.} $$
This proves 3.7.\QED

Now to complete the proof of Theorem 3.1, we must find which affine maps
\[L(\tau)=a\,\tau+b\] give rise to well defined maps of the quotient surface
\[\T/G_n\]. Let \[\omega\] be a primitive \[n$-th root of unity, so that the
rotation \[g(t)=\omega t\] generates \[G_n\]. Then evidently
the points \[L(t)=at+b\] and \[L(g(t))=a\,\omega\, t+b\] represent the same
element of \[\T/G_n\] if and only if
$$  a\,\omega\, t+b~\equiv~ \omega^k(at+b)~~\mod\Lambda\qquad
{\rm for~some~power}\qquad\omega^k~. $$
If this equation is true for some generic choice of \[t\], then it will be true
identically for all \[t\]. Now differentiating with respect to \[t\] we see
that \[k=1\], and substituting \[t=0\] we see that
\[b\equiv \omega\, b~~\mod\Lambda\]. It follows easily that
\[g\circ L=L\circ g\]. Conversely, whenever \[g\] and \[L\] commute it follows
immediately that \[L/G_n\] is well defined. This completes the proof of 3.1.
 \QED\ss

The analogous statement for Chebyshev maps and power maps is the following.

{\QP{\bf Lemma 3.8.} \it If \[f=\Theta\circ L\circ\Theta^{-1}\] is a finite
quotient of an affine map on a cylinder \[\cal C\], then \[f\]
is holomorphically conjugate either to a power map \[z\mapsto z^a\]
or to a Chebyshev map \[\pm\Ch_d\].\ss}

The proof is completely analogous to the proof of 3.1. In fact any such \[f\]
is conjugate to a map of the form \[L/G_n:{\cal C}/G_n\to{\cal C}/G_n\],
where \[L\] is an affine map of the cylinder \[\cal C\] and
where \[n\] is either one (for the power map case) or two (for the Chebyshev
case). Details will be left to the reader.\QED\ss

The following helps to demonstrate the extremely restricted dynamics associated
with finite quotients of affine maps. Presumably
nothing like it is true for any other rational map.\ss

{\QP{\bf Corollary 3.9.} 
\it Let \[f=\Theta\circ L\circ\Theta^{-1}\] be a finite quotient of an affine
map \[L\] which has derivative \[L'=a\]. If \[z\in\wC\ssm{\cal E}_f\] is
a periodic point
with period \[p\ge 1\] and ramification index \[r\ge 1\], then the
multiplier of \[f^{\circ p}\] at \[z\]
is a number of the form \[\mu=(\omega a^p)^r\] where \[\omega^n=1\].\ss}

\ni (The ramification index is described in 3.6 and also in \S4.)
For example for a periodic orbit of ramification \[r=n\] the
multiplier is simply \[a^{p\,n}\]. In the case of a generic periodic orbit
with \[r=1\], the multiplier has the form \[\omega a^p\]. In all cases,
the absolute value \[|\mu|\] is equal to \[|a|^{pr}\].\ss

{\bf Proof of 3.9.} First consider a fixed point \[z_0=f(z_0)\] and let
\[\Theta(\tau_0)=z_0\]. As in 3.6,
we can take \[\zeta=(\tau-\tau_0)^r\] as local uniformizing parameter near
\[z_0\]. On the other hand, since
\[z_0=f(z_0)\] we have \[\tau_0\sim L(\tau_0)\], or in other words \[\tau_0=
\omega\,L(\tau_0)\] for some \[\omega\in G_n\]. Thus \[f\] lifts to the
linear map
$$ \tau-\tau_0~\mapsto~\omega L(\tau)-\omega L(\tau_0)\=
 \omega\,a\,(\tau-\tau_0)~. $$
Therefore, in terms of the local coordinate \[\zeta\] near \[z_0\], we have
the linear map
\[\zeta\mapsto (\omega\,a)^r\zeta\], with derivative \[\mu=(\omega\,a)^r\].
Applying the same argument to the \[p$-th iterates of \[f\] and \[L\], we get a
corresponding identity for a period \[p\] orbit.\QED\bs

{\bf\S4. Flat Orbifold Metrics.} We can give another
characterization of finite quotients of affine maps as follows.\ss

{\bf Definition.} By a {\bit flat orbifold metric\/} on \[\wC\ssm{\cal E}_f\]
will be meant a complete metric which is smooth, conformal,
and locally isometric to
the standard flat metric on \[\C\], except at finitely many ``cone points'',
where it has cone angle of the form \[2\pi/r\].
Here a {\bit cone point\/} with {\bit cone angle\/} \[0<\alpha<2\pi\], is an
isolated singular point of the metric which can be visualized by cutting
an angle of \[\alpha\] out of a sheet of paper and
then gluing the two edges together. (A more formal definition will be left
to the reader.) In the special case where \[\alpha\]
is an angle of the form \[2\pi/r\], we can identify such a cone with the
quotient space \[\C/G_r\] where \[G_r\] is the group of \[r$-th roots of unity
acting by multiplication on the complex numbers, and where the flat metric on
\[\C\] corresponds to a flat metric on the quotient, except at the
cone point.
\midinsert
\FIG{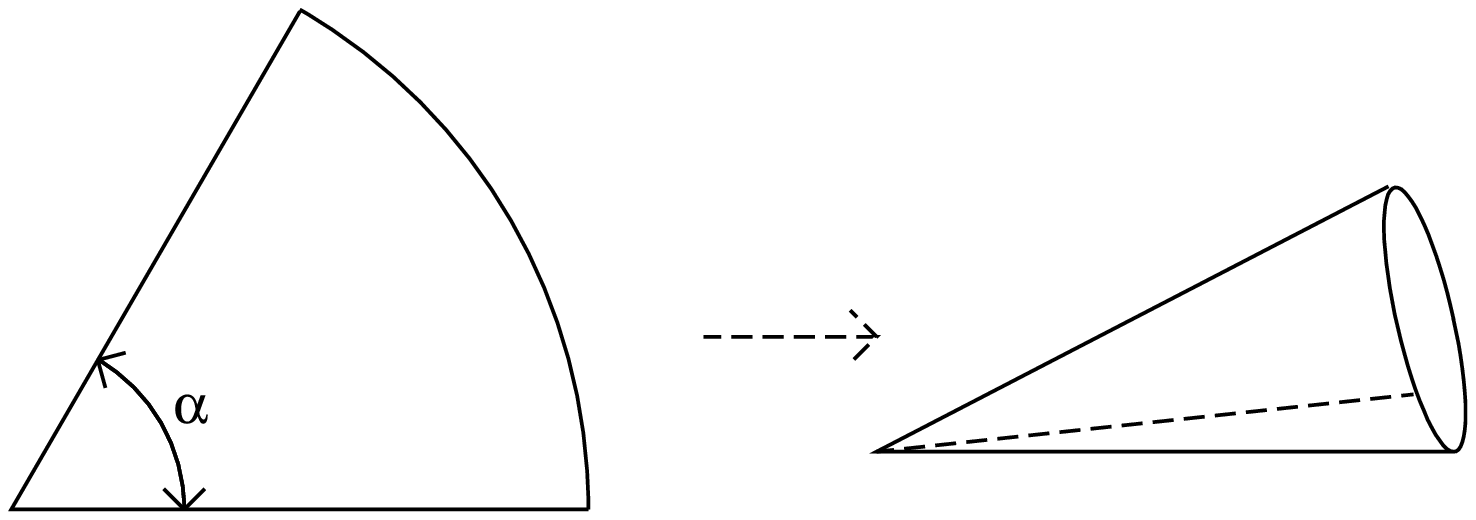}{1.2in}\ss
\cl{\bit Figure 1. Model for a cone point with cone angle \[\alpha\].}
\endinsert

Evidently the canonical flat metric on a torus \[\T\] or cylinder \[\cal C\]
induces a corresponding
flat orbifold metric on the quotient \[\T/G_n\] of Theorem 3.1 or
the quotient \[{\cal C}/G_n\] of Lemma 3.8.
Thus near any non-cone point we can
choose a local coordinate \[t\] so that the metric takes the form
\[|\d t|^2\].
I will say that such a metric {\bit linearizes\/} the map \[f\] since, in terms
of such preferred local coordinates, \[f\] is an affine map with constant
derivative.\footnote{$^4$}{\tp
The more usual terminology for a map \[f\]
which is linearized by a flat orbifold metric would be that \[f\] has
``{\bit parabolic\/}'' or ``{\bit Euclidean{\rm ''} orbifold\/}. Following
Thurston, for
any postcritically finite \[f\] there is a smallest function \[r(z)\ge 1\]
on \[\wC\ssm{\cal E}_f\] such that \[r(f(z))\] is a multiple of
\[d_f(z)\,r(z)\] for every \[z\]. Furthermore there is an essentially unique
complete orbifold metric of constant curvature \[\le 0\] on
\[\wC\ssm{\cal E}_f\] with \[r(z)\]
as ramification index function. The curvature is
zero if and only if (7) is satisfied. (See [DH] or [M3].)}
(An equivalent property is that \[f\] maps any curve of length
\[\delta\] to a curve of length \[k\,\delta\] where \[k=|a|>1\] is
constant.) A converse statement is also true:

{\QP{\bf Theorem 4.1.} \it A rational map \[f\] is a finite quotient of an
affine map if and only if it is linearized by some flat orbifold
metric, or if and only if there exists an integer
valued ``ramification index'' function \[r(z)\] on \[\wC\ssm{\cal E}_f\]
satisfying the identity
$$	r(f(z))\=d_f(z)\;r(z)\qquad{\it for~ all} \quad z~~,\eqno(7) $$
with \[r(z)=1\] outside of the postcritical set of \[f\].\ms}

{\bf Proof in the Latt\`es case.} First suppose that \[f\] is a finite
quotient of an affine map on a torus, conformally conjugate to the
quotient map \[L/G_n:\T/G_n\to\T/G_n\].
If \[\tau_0\] is a critical point of the projection
\[\Theta:\T\to\T/G_n\cong
\widehat \C\] with local degree \[d_\Theta(\tau_0)=r\], then the subgroup
consisting of elements of \[G_n\] which fix \[\tau_0\] must be generated by
a rotation through \[2\pi/r\] about \[\tau_0\]. Hence the flat metric on \[\T\]
pushes forward to a flat metric on \[\T/G_n\] with \[\tau_0\] corresponding to
a cone point \[z_0=\Theta(\tau_0)\] of angle \[2\pi/r\]. This
integer \[r=r(z_0)>1\] is called the ramification index of the cone
point. Setting \[r(z)=1\] if \[z\] is not a cone point, we see that
\[r(\theta(\tau))\] can be identified with the local degree
\[d_\Theta(\tau)\] in all cases.
There may be several different points in \[\Theta^{-1}(z)\], but \[\Theta\]
must have the same local degree at all of these points, since the angle
at a cone point is uniquely defined, or by 3.9. With these notations, the
required equation (7) is just a restatement of equation (3) of \S3.

Conversely, suppose that (7) is satisfied. {\it It follows from
this equation that \[f\] is postcritically finite.\/} In fact we can express
the postcritical set \[P_f\] as a union
\[P_1\cup P_2\cup\cdots\] of disjoint finite subsets, where $$P_1\=f(C_f)\qquad
{\rm and}\qquad P_{m+1}\=f(P_m)\ssm(P_1\cup\cdots\cup P_m)~.$$
Let \[|P_m|\] be the number
of elements in \[P_m\]. Since \[f(P_m)\supset P_{m+1}\], the sequence\break
\[|C_f|\,\ge\, |P_1|\,\ge\, |P_2|\,\ge\,|P_3|\,\ge\,\cdots~\]
must eventually stabilize. Therefore we can choose an integer \[m\] so that
\[P_k\] maps bijectively onto \[P_{k+1}\] for \[k\ge m\], and must prove that
the number of elements \[|P_m|=|P_{m+1}|=\cdots\,\] is zero.
Note that each point of \[P_{m+1}\] has \[d_f\]
distinct preimages, where \[d_f\ge 2\] by our standing hypothesis. Thus if
\[|P_{m+1}|\ne 0\] there would
exist some point \[z\not\in P_m\] with \[f(z)\in P_{m+1}\]. In fact it would
follow
that \[z\not\in C_f\cup P_f\]. For if \[z\] were in \[C_f\cup P_1\cup\cdots
\cup P_{m-1}\] then \[f(z)\] would be in some \[P_k\] with \[k< m+1\], while
if \[z\] were in \[P_k\] with \[k>m\] then \[f(z)\] would be in \[P_{k+1}\]
with \[k+1>m+1\], contradicting the hypothesis that \[f(z)\in P_{m+1}\]
in either case. The
existence of a point \[z\not\in C_f\cup P_f\] with \[f(z)\in P_f\]
clearly contradicts equation (7).

Let me use the notation \[M\] for the Riemann sphere \[\wC\] together with
the ``orbifold structure'' determined by the function
\[r:\wC\to\{1,2,3,\ldots\}\]. The {\bit universal covering orbifold\/}
\[\widetilde M\] can be characterized
as a simply connected Riemann surface together with a holomorphic branched
covering map \[\theta:\widetilde M\to M=\wC\] such that, for each
\[\widetilde z\in \widetilde M\], the local degree \[d_\theta(\widetilde z)\]
is equal to the prescribed ramification index \[r(\theta(\widetilde z))\].
Such a universal covering associated with a function \[r:M\to\{1,2,3,\ldots\}\]
exists whenever the number of \[z\] with \[r(z)>1\] is finite with at least
three elements. (See for example [M3, Lemma E.1].)

For this proof only, it will be
convenient to choose some fixed point of \[f\] as base point \[z_0\in M\].
Using equation (7), there is no obstruction to lifting \[f\] to a
holomorphic map \[\widetilde f\] which maps the Riemann surface
\[\widetilde M\] diffeomorphically into itself, with
no critical points. Furthermore, we can choose
\[\widetilde f\] to fix some base point \[\widetilde z_0\] lying over \[z_0\].
The covering manifold \[\widetilde M\] cannot be a compact surface, necessarily
of genus zero, since then \[\widetilde f\] and hence \[f\] would have degree
one, contrary to the standing hypothesis that \[d_f\ge 2\]. Furthermore,
since the base point in \[\widetilde M\] is
strictly repelling under \[\widetilde f\], it follows that \[\widetilde M\]
cannot be a hyperbolic surface.
Therefore \[\widetilde M\] must be conformally isomorphic to the complex
numbers \[\C\], and \[\widetilde f\] must correspond to a linear map \[L\]
from \[\C\] to itself. Evidently the standard flat metric on
\[\C\cong \widetilde M\] now gives rise to a flat orbifold
metric on \[M=\wC\] which linearizes the map \[f\].

Finally, suppose that we start with a flat orbifold metric on \[\wC\]
which linearizes the rational map \[f\]. The preceding discussion shows
that \[f\] lifts to a linear map \[\widetilde f\] on the universal covering
orbifold \[\widetilde M\]. Let \[\widetilde G\] be the group of {\bit deck
transformations\/} of \[\widetilde M\], that is homeomorphisms \[\widetilde g\]
from \[\widetilde M\] to itself which cover the identity map of \[M\], so
that \[\theta\circ\widetilde g=\theta\]. Then the quotient surface
\[\widetilde M/\widetilde G\] can be identified with \[M=\wC\].
If \[\Lambda\subset\widetilde G\]
is the normal subgroup consisting of those deck transformations which are
translations of \[\widetilde M\cong\C\], then
the quotient group \[G=\widetilde G/\Lambda\] is a finite group of
rotations with order equal to the least common multiple of the ramification
indices. It follows that the quotient \[\T=\widetilde M/\Lambda\] is a torus,
and hence that
\[f\] is a finite quotient of an affine map of this torus.\ss

The proof of 4.1 in the Chebyshev and power map cases is similar and will be
omitted.\QED\ms

{\bf Remark 4.2.} Note that the construction of the torus \[\T\], the group
\[G_n\] and the affine map \[L\] from the rational map \[f\] satisfying (7)
is completely canonical, except for the choice of lifting for \[f\].
For example, when there are four postcritical points, the conformal
conjugacy class of the torus \[\T\] is completely determined by the set
of postcritical points, and in fact by the cross-ratio of these four points.

However, to make an explicit classification we must note the following.

\[\bullet\] We want to identify the torus \[\T\] with some quotient
\[\C/\Lambda\]. Here, \[\Lambda\] is unique only up to multiplication
by a non-zero constant; but we can make an explicit and
unique choice by taking
\[\Lambda\] to be the lattice \[\Z\oplus\gamma\Z\] spanned by \[1\]
and \[\gamma\],  where \[\gamma\] belongs to the Siegel region
$$\eqalign{& |\gamma|\ge 1~,\qquad|\Re(\gamma)|\le 1/2~,\qquad
 \Im(\gamma)>0~,\cr &{\rm with}\quad
\Re(\gamma)\ge 0\quad {\rm whenever}\quad|\gamma|=1\quad {\rm or}
\quad|\Re(\gamma)|=1/2~.} \eqno(8)$$
With these conditions,
\[\gamma\] is uniquely determined by the conformal isomorphism class of \[\T\].
We will describe the
corresponding \[\Lambda=\Z\oplus\gamma\Z\] as a {\bit normalized lattice\/}.

\[\bullet\] For specified \[\Lambda\], we still need to make some choice
of conformal isomorphism\break \[\upsilon:\C/\Lambda\to\T\]. In most cases,
\[\upsilon\] depends only on a choice of base point  \[\upsilon(0)\in\T\],
up to sign. However, in the special case where \[\T\] admits a \[G_3\] (or
\[G_4\]) action, we can also multiply \[\upsilon\] by a cube (or fourth)
root of unity. As in \S3,
it will be convenient to choose one of
the fixed points of the \[G_n\] action as a base point in \[\T\].

\[\bullet\] The lifting \[L(t)=at+b\] of the map \[f\] to the torus is well
defined only up to the action of \[G_n\]. In particular, we are always free
to multiply the coefficients \[a\] by an \[n$-th root of unity.

We will deal with all of these ambiguities in \S5.\ms

Here is an interesting consequence of 4.1. Let \[f\] and \[g\] be rational
maps.

{\QP{\bf Corollary 4.3.} \it Suppose that there is a holomorphic semiconjugacy
from \[f\] to \[g\], that is, a 
non-constant rational map \[h\] with \[h\circ f=g\circ h\].
Then \[f\] is a finite quotient of an affine map if and only if \[g\]
is a finite quotient of an affine map.\ss}

{\bf Proof.} It is not hard to see that \[h^{-1}({\cal E}_g)={\cal E}_f\],
so that \[h\] induces a proper map from \[\wC\ssm{\cal E}_f\] to
\[\wC\ssm{\cal E}_g\]. Now if \[f\] is a finite quotient of an affine map
\[L\], say \[f=\Theta\circ L\circ \Theta^{-1}\], then it follows immediately
that
\[g=(h\circ\Theta)\circ L\circ(h\circ\Theta)^{-1}\]. Conversely, if \[g\] is
such a finite quotient, then there is
a flat orbifold structure on \[\wC\ssm{\cal E}_g\] which linearizes \[g\],
and we can lift easily to a flat orbifold structure on \[\wC\ssm{\cal E}_f\]
which linearizes \[f\].\QED\ms

In order to 
classify all possible flat orbifold structures on the Riemann sphere,
we can use a piecewise linear form of  the Gauss-Bonnet Theorem. For this
lemma only, we allow cone angles which are greater than \[2\pi\].

{\QP{\bf Lemma 4.4.} \it If a flat
metric with finitely many cone points on a compact Riemann surface \[S\]
has cone angles \[\alpha_1\,,\,\ldots\,,\,\alpha_k\], then
$$	(2\pi-\alpha_1)\,+\,\cdots\,+\,(2\pi-\alpha_k)\= 2\pi\,\chi(S)~,
 \eqno(9) $$
where \[\chi(S)\] is the Euler characteristic. In particular, if \[\alpha_j=
2\pi/r_j\] and if \[S\] is the Riemann sphere with \[\chi(S)=2\],
then it follows that \[~\sum\,(1-1/r_j)\=2\].\ss}

{\bf Proof.} Choose a rectilinear triangulation, where the cone points will
necessarily be among the vertices. Let \[V\] be the number of vertices, \[E\]
the number of edges, and \[F\] the number of faces (i.e., triangles). Then
\[2E=3F\] since each edge bounds two triangles and each triangle has three
edges. Thus
$$	\chi(S)\=V-E+F\=V-F/2~.\eqno(10)$$
The sum of the internal angles of all of the
triangles is clearly equal to \[\pi F\]. On the other hand, the \[j$-th cone
point contributes \[\alpha_j\] to the total, while each non-cone vertex
contributes \[2\pi\]. Thus
$$ \pi F\=\alpha_1+\cdots+\alpha_k+2\pi( V-k)~. \eqno(11)$$
Multiplying equation (10) by \[2\pi\] and using
(11), we obtain the required equation (9).\QED\ss

{\QP{\bf Corollary 4.5.} \it
The collection of ramification indices for a flat orbifold metric on the
Riemann sphere
must be either \[\{2,2,2,2\}\] or \[\{3,3,3\}\] or \[\{2,4,4\}\] or
\[\{2,3,6\}\].
In particular, the number of cone points
must be either four or three.\ss}

{\bf Proof.} Using the inequality \[1/2\le (1-1/r_j)<1\],
it is easy to check that the required equation
$$	\sum_j \;(1-1/r_j) \=\chi(\widehat\C)\=2~, $$
has only these solutions in integers \[r_j>1\].\QED\ss

{\bf Remark 4.6.} If \[z\in\wC\] corresponds to a fixed point for
the action of the group \[G_n\] on the torus, then the ramification index
\[r(z)\] is evidently
equal to \[n\]. For any other point, it is some divisor of \[n\]. Thus
the order \[n\] of the rotation group \[G_n\] can be identified with the
least common multiple (or the maximum) of the various ramification indices
as listed in 4.5.\ss

{\bf Remark 4.7.} To deal with the case of a map \[f\] which has exceptional
points, we can
assign the ramification index \[r(z)=\infty\] to any exceptional point
\[z\in{\cal E}_f\]. If we allow such points, then the equation
\[\sum(1-1/r_j)=2\] has two further solutions, namely:
\[\{\infty\,,\,\infty\}\]
corresponding to the power map case, and \[\{2\,,\,2\,,\,\infty\}\]
corresponding to the Chebyshev case.\ss

Combining Corollary 4.5 with equation (7), we get an easy
characterization of Latt\`es maps in two of the four cases.

{\QP{\bf Corollary 4.8.} \it A rational map with four postcritical points is
Latt\`es if and only if every critical point is simple $($with local degree
two$)$
and no critical point is postcritical. Similarly, a rational map with three
postcritical points is Latt\`es of type \[\{3,3,3\}\] if and only if
every critical point has local degree three and none is postcritical.\ss}

The proof is easily supplied.\QED\ms

\pageinsert
\FIG{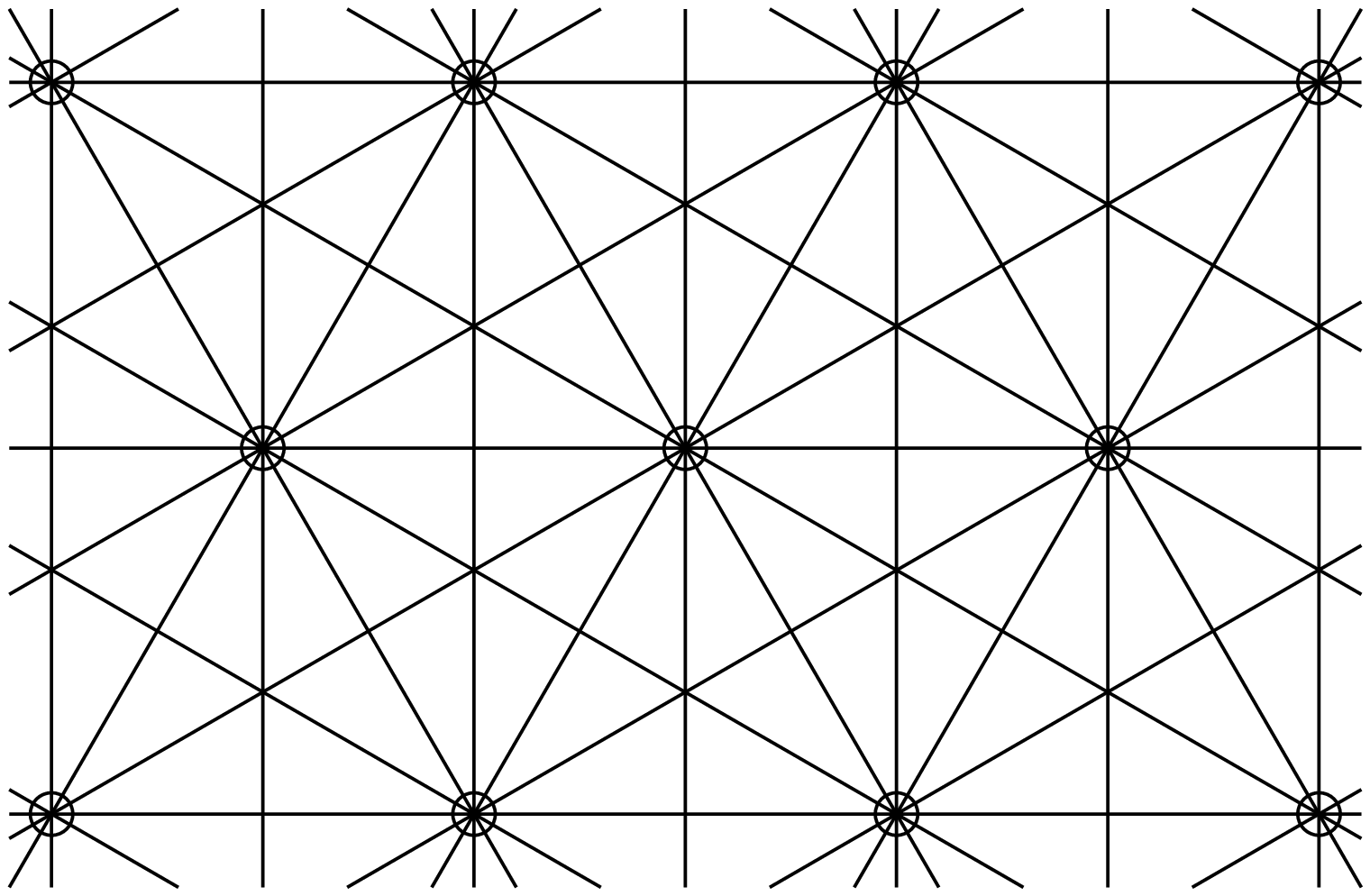}{2.4in}
{\QP\bit Figure 2. The \[\{2,3,6\}$-tiling of the plane. In each of these
diagrams, the points  of ramification \[n\] have been marked,
with circles around the lattice points.\ms}

\FIG{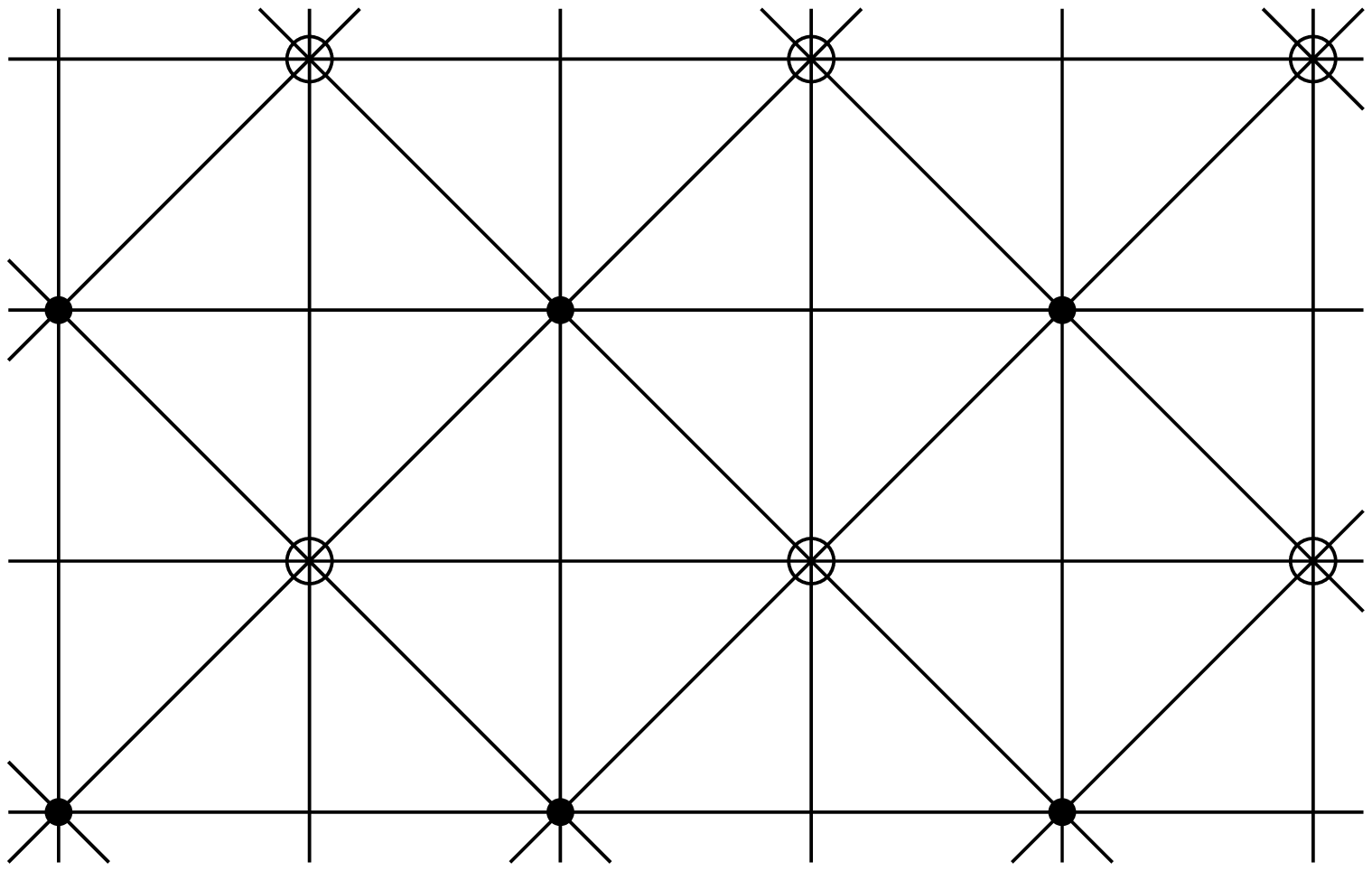}{2.4in}
\cl{\bit Figure 3. The \[\{2,4,4\}$-tiling.}\ms

\FIG{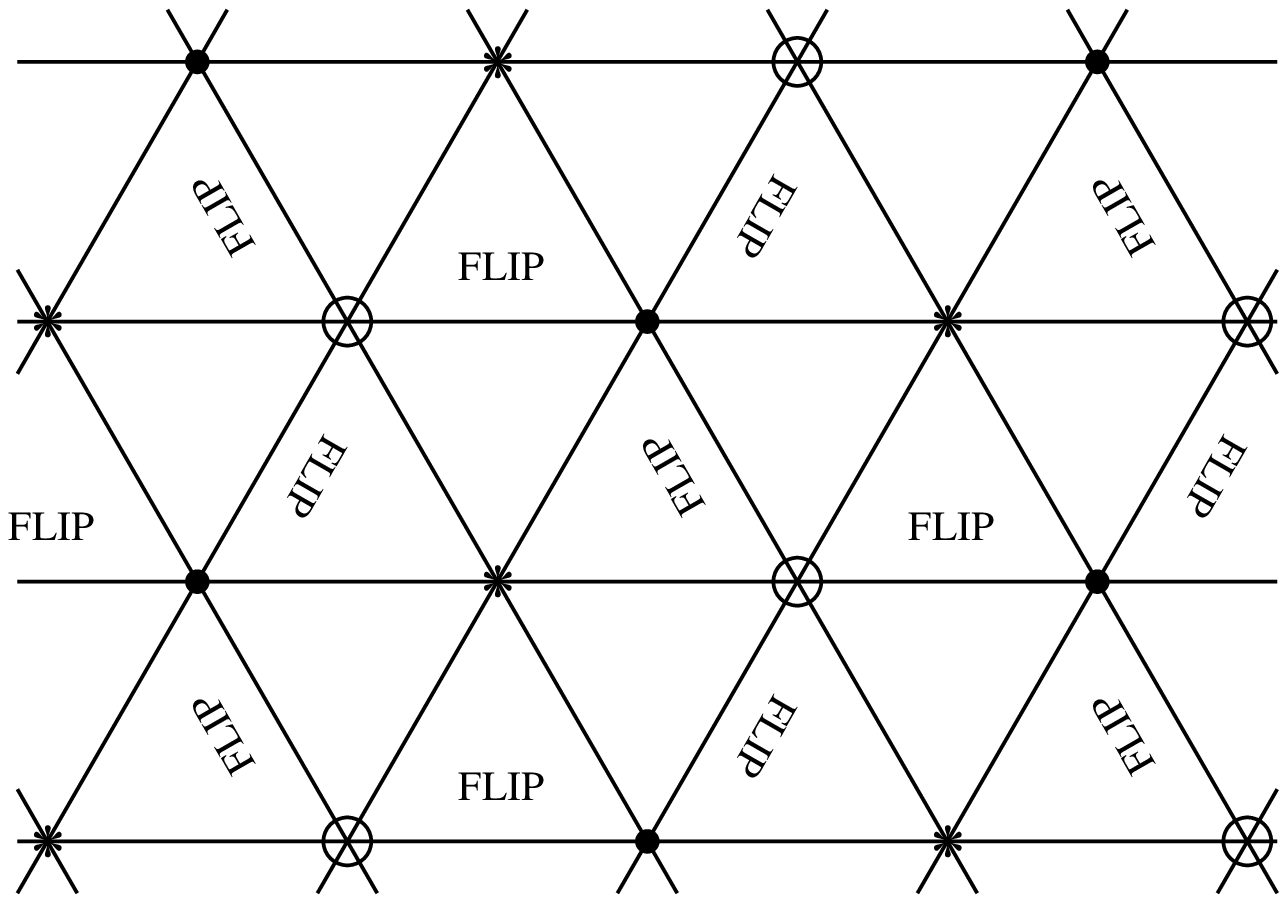}{2.7in}
\cl{\bit Figure 4. The
\[\{3,3,3\}$-tiling, with one tile and its images under \[\widetilde G_3\]
labeled.}
\endinsert

We conclude this section with a more precise description of the possible
crystallographic groups \[\widetilde G\] acting on \[\C\], and
of the corresponding  orbifold geometries on \[\C/\widetilde G\cong\T/G_n\].
We first look at the cases \[n\ge 3\] where there are exactly three
cone points in \[\T/G_n\] or equivalently
three postcritical points for any associated Latt\`es map. Thus the
collection of ramification indices must be
either \[\{2,3,6\}\] or \[\{2,4,4\}\] or \[\{3,3,3\}\].
Each of these three possibilities is associated
with a rigidly defined flat orbifold geometry
which can be described as follows. Join each pair of cone points by a minimal
geodesic. Evidently these geodesics cannot cross each other; and no geodesic
can pass through a cone point since our cone angles are strictly less than
\[2\pi\]. In this way, we obtain
three edges which cut our locally flat manifold into
two Euclidean triangles. Since these two triangles have the same edges,
they must be precise mirror images of each other. In particular, the two
edges which meet at a cone point of angle \[2\pi/r_j\] must cut it into
two Euclidean angles of \[\pi/r_j\]. Passing to the branched covering space 
\[\T\] or its universal covering
\[\widetilde \T\], we obtain a tiling of the torus or the Euclidean
plane\footnote{$^5$}{\tp More generally, for any triple of integers
\[r_j\ge 2\] there is an associated tiling, either of the Euclidean
or hyperbolic plane or of the 2-sphere depending on the sign of
\[1/r_1+1/r_2+1/r_3-1\]. See for example [M1].}
by triangles with angles \[\pi/r_1\,,\,\pi/r_2\] and \[\pi/r_3\]. These
tilings are illustrated in Figures 2, 3, 4.

In each case, each pair of adjacent triangles are mirror images of each other,
and together form a fundamental domain for the action of
the group of Euclidean motions \[\widetilde G_n\]
on the plane, or for the action of \[G_n\] on the torus. For each vertex
of this diagram, corresponding to a cone point of angle \[2\pi/r_j\],
there are \[r_j\] lines through the vertex, and hence \[2r_j\] triangles
which meet at the vertex. The subgroup of \[\widetilde G_n\] (or \[G_n$)
which fixes such a
point has order \[r_j\] and is generated by a rotation through the
angle \[\alpha_i=2\pi/r_j\].

The subgroup \[\Lambda\subset\widetilde G_n\] consists of all translations
of the plane which belong to \[\widetilde G_n\].
Recall from 4.4 that the integer \[n\] can be
described as the maximum of the \[r_j\]. The \[2n\] triangles which
meet at any maximally complicated vertex form a fundamental domain
for the action of this subgroup \[\Lambda\]. In the \[\{2,3,6\}\] and
\[\{3,3,3\}\] cases, this fundamental domain is a regular hexagon, while in the
\[\{2,4,4\}\] case it is a square. In all three cases, the torus \[\T\]
can be obtained by identifying opposite faces of this fundamental domain
under the appropriate translations. Thus when \[n\ge 3\] the torus
\[\T\] is uniquely determined by \[n\], up to conformal diffeomorphism.

In the \[\{2,3,6\}\] case, the integers \[r_j\] are all distinct, so it is
easy to distinguish the three kinds of vertices. However, in the
\[\{2,4,4\}\] case there are two different kinds of vertices of index 4.
In order to distinguish them, one kind has been marked with dots and the
other with circles. Similarly in the \[\{3,3,3\}\] case, the three kinds of
vertices have been marked in three different ways.
In this last case, half of the triangles have also been labeled.
In all three cases, the points of the lattice \[\Lambda\], corresponding
to the base point in \[\T\], have been circled.
For all three diagrams, the group \[\widetilde G_n\] can be descried as the
group of all rigid Euclidean motions which carry the marked diagram to itself,
and the lattice  \[\Lambda\] can be identified with the subgroup
consisting of translations which carry this marked diagram to itself.\ss

\midinsert
\FIG{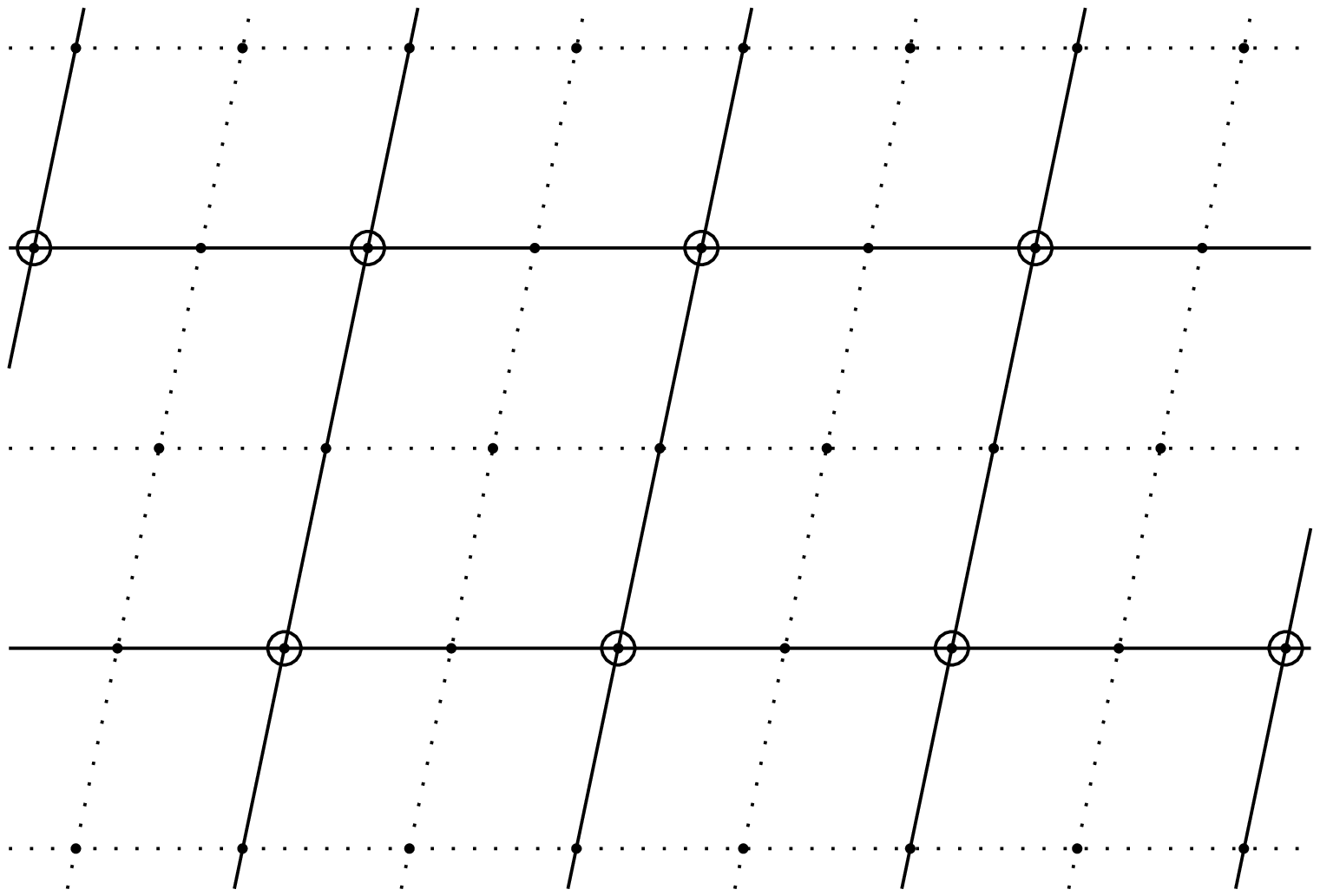}{1.8in}\bs
\cl{\bit Figure 5. A typical \[\{2,2,2,2\}\] tiling of the plane.}
\endinsert

The analogue of Figures 2, 3, 4 for a typical orbifold of type
\[\{2,2,2,2\}\] is a tiling of the plane by parallelograms associated
with a typical lattice \[\Lambda=\Z\oplus\gamma\Z\],
as illustrated in Figure 5. All of the vertices in this figure represent
critical points for the projection map\break
 \[\theta:\C\to\widehat\C\]. Again lattice points have
been circled. Any two adjacent small parallelograms form a fundamental region
for the action of the group \[\widetilde G_2\], which consists of \[180^\circ\]
rotations around the vertices, together with lattice translations. The
four small parallelograms
adjacent to any vertex forms a fundamental domain under lattice translations.

\midinsert
\FIG{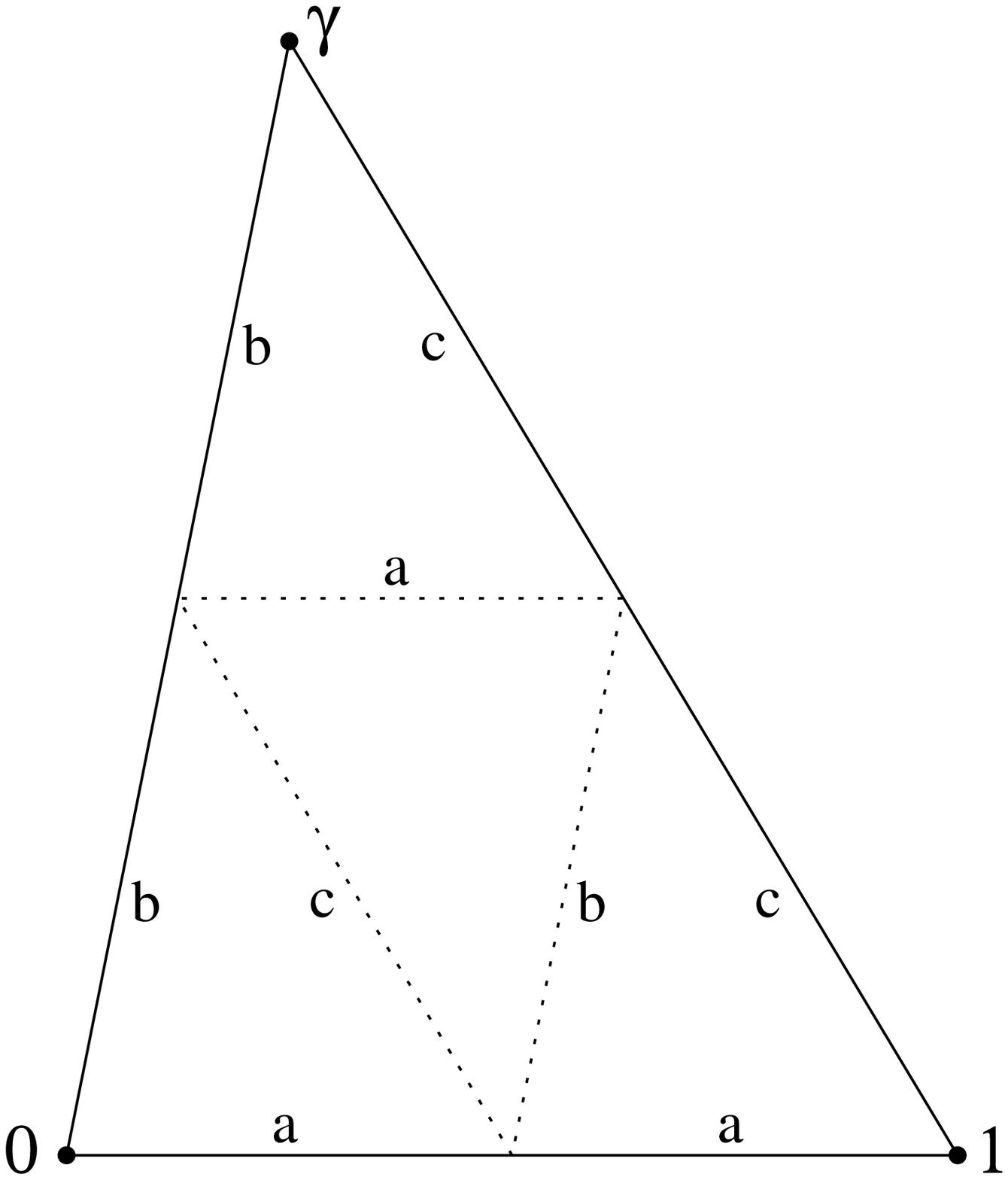}{2.5in}
\cl{\bit Figure 6. Illustrating the orbifold structure of \[\T/G_2\].}
\endinsert

In most cases, the corresponding flat orbifold is isometric to some
tetrahedron in Euclidean space. (Compare [De].) For example, consider the case
where the invariant \[\gamma\] in the Siegel domain (8) satisfies
\[0<\Re(\gamma)<1/2\]. Then the triangle with vertices \[0\,,\,1\] and
\[\gamma\] has all angles acute, and also serves
as a fundamental domain for the
action of the crystallographic group \[\widetilde G_2\] on \[\C\]. Joining
the midpoints of the edges of this triangle, as shown in Figure 6, we can
cut this triangle up into four similar triangles. Now fold along the dotted
lines and bring the three corner triangles up so that the three vertices
\[0\,,\,1\] and \[\gamma\] come together. In this way, we obtain a tetrahedron
which is isometric to the required flat orbifold \[\T/G_2\]. The construction
when \[-1/2\le\Re(\gamma)<0\] is the same, except that we use \[-1\] and \[0\]
in place of \[0\] and \[1\]. The tetrahedrons which can be obtained in this way
are characterized by the property that
opposite edges have equal length, or by the property that there is a
Euclidean motion carrying any vertex to any other vertex. In most cases this
Euclidean motion is uniquely determined; but in the special case
where we start with an
equilateral triangle, with \[\gamma=\omega_6\], we obtain a regular
tetrahedron which has extra symmetries. !!!!

In the case \[\Re(\gamma)=0\]
where \[\gamma\] is pure imaginary, this tetrahedron
degenerates and the orbifold can be described rather as the ``double''
of the rectangle which has vertices
$$0\;,~~1/2\;,~~\gamma/2\;,~~(1+\gamma)/2~.$$
Again, in most cases there is a unique orientation preserving isometry carrying
any vertex to any other vertex; but in the special case of a square,
with \[\gamma=i\], there are extra symmetries.

\bs

{\bf\S5. Classification.} By taking a closer look at
the arguments in sections 3 and 4, we can give a precise classification
of Latt\`es Maps. (Compare [DH, \S9].) It will be convenient to introduce the
notation
$$\omega_n\=\exp(2\pi i/n)~,\eqno(12) $$
for the standard generator of the cyclic group \[G_n\]. Thus
$$ \omega_2=-1\,,\qquad\omega_3=(-1+i\sqrt 3)/2\,,
\qquad\omega_4=i\,,\qquad \omega_6=\omega_3+1~.$$
As usual, we consider a Latt\`es map which is conjugate to
\[L/G_n:\T/G_n\to\T/G_n\], where \[\T\cong\Z/\Lambda\] and where \[L(t)=at+b\].
Here it will be convenient to think of \[b\] as a complex number, well
defined modulo \[\Lambda\].

{\QP{\bf Theorem 5.1.} \it Such a Latt\`es map \[f\] is uniquely determined
up to conformal conjugacy by the following four invariants.

\ni\[\bullet\] First: the integer \[n\], equal to 2, 3, 4, or 6.

\ni\[\bullet\] Second: the complex number \[a^n\], with \[|a|^2\] equal to
the degree of \[f\].

\ref\[\bullet\] Third: the lattice \[\Lambda\], which we may take to have the
form  \[\Lambda=\Z\oplus\gamma\Z\] with \[\gamma\] in the Siegel region $(8)$.
This lattice must satisfy the conditions that
\[\omega_n\Lambda=\Lambda\] and \[a\Lambda\subset\Lambda\].
Let \[k\] be the largest integer such that \[\omega_k\Lambda=\Lambda\].

\ni\[\bullet\] Fourth: the product \[\;(1-\omega_n)\,b\in
\Lambda\] modulo the sublattice $$(1-\omega_n)\Lambda\,+\,(a-1)\Lambda
\;\subset\;\Lambda~,$$ up to multiplication by \[G_k\] with \[k\]
as above. This last
invariant is zero if and only if the map \[f\] admits a fixed point of
maximal ramification index \[r=n\], or equivalently a fixed point of
multiplier \[\mu=a^n\].\ms}

\ni For most lattices we have \[k=2\], so that
the image of  \[(1-\omega_n)\,b\] in the quotient group
$$ \Lambda\Big/\Big((1-\omega_n)\Lambda\,+\,(a-1)\Lambda\Big) \eqno (13)$$
is invariant up to sign.
However, in the special case where \[\Lambda\] has \[G_4\] or \[G_6\] symmetry,
so that \[\gamma=\omega_4\] or \[\gamma=\omega_6$, 
this image is invariant only up to
multiplication by \[\omega_4\] or \[\omega_6\] respectively.

Note that the first invariant \[n\], equal to the greatest common divisor of
the ramification indices, can easily be computed by looking at the orbits of
the critical points of \[f\], using formula (7) of \S4. The invariant \[a^n\]
can be computed from the multiplier \[\mu\] at any fixed point, since the
equation \[\mu=(\omega a)^r\] of Corollary 3.9, with \[\omega^n=1$,
implies that \[\mu^{n/r} =a^n\]. It follows from this equation
that \[\mu=a^n\] if and only if \[r=n\].

The cases with \[n\ge 3\] are somewhat easier than the case \[n=2\].
In fact a normalized lattice \[\Lambda\] with \[G_3\] or \[G_6\] symmetry
is necessarily equal to \[\Z[\omega_3]=\Z[\omega_6]\], and the condition \[a\Lambda\subset
\Lambda\] is satisfied if and only if \[a\in
\Z[\omega_6]\]. Similarly, \[G_4\] symmetry implies that \[\Lambda=\Z[i]\],
and the possible choices for \[a\] are just the elements of \[\Z[i]\],
always subject to the standing requirement that \[|a|>1\].

{\QP{\bf Corollary 5.2.} \it If \[n\ge 3\], then the conformal conjugacy
class of \[f\] is completely determined by the numbers
\[n\] and \[a^n\] where \[a\in\Z[\omega_n]\],
together with the information as to whether
\[f\] does or does not have a fixed point of multiplier \[\mu=a^n\].\ms}

Evidently there is such a fixed point if and only if
\[(1-\omega_n)\,b\] is congruent to zero modulo
\[(1-\omega_n)\Lambda\,+\,(a-1)\Lambda\]. (When \[n=6\] there is
necessarily such a fixed point.)\ss

The proof of 5.1 and 5.2 will be based on the following.\ss

{\QP{\bf Lemma 5.3.} \it The additive subgroup of \[\T=\C/\Lambda\] consisting
of elements fixed by the action of \[G_n\] is canonically isomorphic to the
quotient group \[\Lambda/(1-\omega_n)\Lambda\], of order \[|1-\omega_n|^2=
4\sin^2(\pi/n)\]. In fact, the correspondence \[t\mapsto(1-\omega_n)\,t\]
maps the group of torus elements fixed by \[G_n\] isomorphically onto
this quotient group.\ss}

{\bf Proof.} The required equation \[\omega_n t\equiv t~~\mod\Lambda\] is
equivalent to \[(1-\omega_n)\,t\in\Lambda\], and the conclusion follows
easily.\QED

Note that points of \[\T\] fixed by the action of \[G_n\] correspond to
points in the quotient sphere \[\T/G_n\] of maximal ramification index \[r=n\].
As a check,  in the four cases \[\{2,2,2,2\}\,,~~\{3,3,3\}\,,~~\{2,4,4\}\],
and \[\{2,3,6\}\], there are clearly 4, 3, 2, and 1 such points respectively.
This number is equal to \[4\sin^2(\pi/n)\] in each case.\ss

{\bf Proof of 5.1.} It is clear that the numbers \[n\,,~\gamma\,,~a^n\],
and \[b\] completely determine the map \[L/G_n:\T/G_n\to\T/G_n\]. In fact
\[\gamma\] determines the torus \[\T\], and the power
\[a^n\] determines \[a\] up to multiplication by \[n$-th roots of unity.
But we can multiply \[L\] and hence \[a\]
by any \[n$-th root of unity without changing the quotient \[\L/G_n\].
Since the numbers \[n\,,~a^n\], and \[\gamma\] are uniquely determined by
\[f\] (compare the discussion above), we need only study
the extent to which \[b\] is determined by \[f\].

Recall from Theorem 3.1 that the map \[L(t)=at+b\] must commute with
multiplication by \[\omega_n\]. That is
$$ L(\omega_nt)~\equiv\omega_n L(t)\qquad\mod\Lambda~.$$
Taking \[t=0\] it follows that \[b\equiv\omega_n b~~\mod\Lambda\], or in
other words \[(1-\omega_n)\,b\in\Lambda\], as required.

Next recall that we are free to choose any fixed point \[t_0\] for the action
of \[G_n\] on \[\T\] as base point. Changing the base point to \[t_0\] in place
of \[0\] amounts to replacing \[L(t)\] by the conjugate map
\[L(t+t_0)-t_0 \= at\,+\,b'\] where
$$ b'\=b+(a-1)\,t_0~,\quad{\rm and~ hence} \quad
(1-\omega_n)\,b'=(1-\omega_n)\,b+(a-1)(1-\omega_n)\,t_0~.$$
Since the product \[(1-\omega_n)\,t_0\] can be a completely arbitrary element
of \[\Lambda\], this means that we can add a completely arbitrary
element of \[(a-1)\Lambda\] to the product \[(1-\omega_n)\,b\] by a
change of base point. Thus the residue class
$$   (1-\omega_n)\,b~\in~
 \Lambda\Big/\Big((1-\omega_n)\Lambda\,+\,(a-1)\Lambda\Big)~,  $$
together with \[n\,,\,\gamma\], and \[a^n\], suffices to determine the
conjugacy class of \[f\]. 
However, we have not yet shown that this residue class
is an invariant of \[f\], since we must also consider automorphisms of \[\T\]
which fix the base point. Let \[\omega\] be any root of unity which satisfies
\[\omega\Lambda=\Lambda\]. Then \[L(t)=at+b\] is conjugate to the map
\[\omega L(t/\omega)=at+\omega b\]. In most cases, we can only choose \[\omega
=\pm 1\]. (The fact that we are free to change the sign of \[b\]
is irrelevant when \[n\] is even, but will be important in the case \[n=3\].)
However, if \[\Lambda=\Z[\omega_6]\] then we can choose \[\omega\] to be
any power of \[\omega_6\], and
if \[\Lambda=\Z[i]\] then we can choose \[\omega\] to be
any power of \[i\]. Further details of the proof are straightforward.\QED\ss

{\bf Proof of 5.2.} In the cases \[n\ge 3\], we have noted that \[\Lambda\]
is necessarily equal to \[\Z[\omega_n]\].
Furthermore, for \[n=3\,,\,4\,,\,6\], the quotient group \[\Lambda/
(1-\omega_n)\Lambda\] is cyclic of order \[3\,,\,2\,,\,1\] respectively.
Thus this group has at most one non-zero element, up to sign. The conclusion
follows easily.\QED\ms

The discussion of Latt\`es maps of type \[\{2,2,2,2\}\] will be divided into
two cases according as \[a\in\Z\] or \[a\not\in\Z\]. First suppose that
\[a\not\in\Z\].\ss

{\bf Definition.}
A complex number \[a\] will be called an {\bit imaginary quadratic integer\/}
if it satisfies an equation
\[\;	a^2\,-\,q\,a\,+\,d=0 \;\]
with integer coefficients and with \[q^2<4d\], so that
$$a\=\Big(q\pm\sqrt{q^2-4d}\,\Big)/2 \eqno (14)$$
is not a real number. Here \[|a|^2=d\] is the associated degree.
Evidently the imaginary quadratic integers form a discrete subset of the
complex plane. In fact for each choice of \[|a|^2=d\] there are roughly
\[4\sqrt d\] possible choices for \[q\], and twice that number for \[a\].\ss

{\QP{\bf Lemma 5.4.} \it A complex number \[a\] can occur
as the derivative \[a=L'\] associated with an affine torus map if and only
if it is either a rational integer \[a\in\Z\] or
an imaginary quadratic integer. If \[a\in\Z\] then any torus can occur,
but if \[a\not\in\Z\] then
there are only finitely many possible tori up to conformal
diffeomorphism. Furthermore, there is a one-to-one correspondence between
conformal diffeomorphism classes of such tori and
ideal classes in the ring \[\Z[a]\].\ms}

{\bf Proof.} Let \[\T=\C/\Lambda\].
The condition that \[a\Lambda\subset\Lambda\] means that \[\Lambda\] must
be a module over the ring \[Z[a]\] generated by \[a\].
We first show that \[a\] must be an algebraic integer.
Without loss of generality, we may assume that \[\Lambda=\Z\oplus\gamma\Z\]
is a normalized lattice, satisfying the Siegel conditions (8). Thus
\[1\in\Lambda\] and it follows that
all powers of \[a\] belong to \[\Lambda\]. If
\[\Lambda_k\] is the sublattice spanned by \[1\,,\,a\,,\,a^2\,,\,\ldots\,,\,
a^k\], then the lattices \[\Lambda_1\subset\Lambda_2\subset\cdots\subset
\Lambda\]
cannot all be distinct. Hence some power \[a^k\]
must belong to \[\Lambda_{k-1}\], which proves that \[a\] satisfies
a monic equation with integer coefficients, and hence is an algebraic integer.
On the other hand, \[a\] belongs to a quadratic number field since the three
numbers \[1\,,\,a\,,\,a^2\in\Lambda\] must satisfy a linear relation with
integer coefficients.
Using the fact that the integer polynomials form a unique factorization domain,
it follows that \[a\] satisfies a monic degree two polynomial.

Now given \[a\not\in\Z\] we must ask which normalized lattices \[\Lambda\]
are possible. Since
\[a\in\Lambda\], we can write \[a=r+s\,\gamma\] with \[r\,,\,s\in\Z\].
Changing the sign of \[a\] if necessary, we
may assume that \[s>0\]. Taking real and imaginary parts, it follows that
$$	r=\Re(a)\,-\,s\,\Re(\gamma)\qquad{\rm and}\qquad
	s\=\Im(a)/\Im(\gamma) ~.$$
On the other hand, it follows easily from the Siegel inequalities
$$   |\gamma|\ge 1~,\qquad|\Re(\gamma)|\le 1/2~,\qquad \Im(\gamma)>0 $$
that \[\Im(\gamma)\ge\sqrt 3/2\]. Since
\[a\] has been specified, this inequality
yields an upper bound of \[2\,\Im(a)/\sqrt 3\]
 for \[s\], and the inequality
\[|\Re(\gamma)|\le 1/2\] then yields an upper bound
for \[|r|\]. Thus there are only finitely many possibilities for
\[\gamma=(a-r)/s\].

Next note that the product lattice \[{\cal I}=s\,\Lambda
=s\Z\oplus(a-r)\Z\] is contained in the ring
\[\Z[a]\], and is an ideal in this ring since \[a\,{\cal I}\subset {\cal I}\].
Clearly the torus \[\C/\Lambda\] is isomorphic to \[\C/{\cal I}\]. If
\[{\cal I}'\] is another ideal in \[\Z[a]\], note that
\[\C/{\cal I}\cong\C/{\cal I}'\] if and only if \[{\cal I}'=c\,{\cal I}\]
for some constant \[c\ne 0\]. Such a constant must belong to
the quotient field \[\Q[a]\], so by definition this means that \[{\cal I}\] and
\[{\cal I}'\] represent the same ideal class.\QED\ss

For further discussion of maps with
imaginary quadratic \[a\] see 7.2, 8.1 and 8.2 below.
We next discuss the case \[a\in\Z\].

{\bf Definition:} A Latt\`es map
$$L/G_n:\T/G_n\to\T/G_n\qquad{\rm with}\qquad \T=\C/\Lambda $$
will be called {\bit flexible~} if we can vary \[\Lambda\] and \[L\]
continuously so as to obtain other Latt\`es maps which are not conformally
conjugate to it.

{\QP{\bf Lemma 5.5.} \it A Latt\`es map \[L/G_n:\T/G_n\to\T/G_n\] is
flexible if and only if \[n=2\], and the affine map \[L(\tau)=a\tau+b\] has
integer derivative, \[L'=a\in\Z\].\ss}

{\bf Proof.} This follows easily from 5.1 and 5.4.\QED\ss

We can easily classify such maps into a single
connected family provided that the
degree \[a^2\] is even, and into two connected families when \[a^2\] is odd,
as follows. In each case, the coefficients \[a\] and \[b\] will
remain constant but \[\T\] will vary through all possible conformal
diffeomorphism classes.\ss

\[\bullet\] {\bf Maps with postcritical fixed point.}  Let \[\H\] be the upper
half-plane, For each integer \[a\ge 2\] there is a connected family of flexible
Latt\`es maps of degree \[a^2\] parametrized by the half-cylinder \[\H/\Z\],
as follows. Let \[\T(\gamma)\] be the torus
\[\C/(\Z\oplus\gamma\Z)\] where \[\gamma\] varies
over \[\H/\Z\], and let \[L:\T(\gamma)\to\T(\gamma)\] be the map \[L(t)=at\].
Then
$$L/G_2: \T(\gamma)/G_2\to\T(\gamma)/G_2 $$
is the required smooth family of maps, with the image of \[0\in\T\]
as ramified fixed point. If we restrict \[\gamma\] to the Siegel
region (8), then we get a set of representative maps which are unique
up to holomorphic conjugacy.\ss

\[\bullet\] {\bf  Maps
without postcritical fixed point.} The construction in this case
is identical, except that we take \[L(t)=at+1/2\]. If \[a\] is even, this
construction yields nothing new. In fact, 
the quotient group (13) of 5.1 is then trivial, and it follows that
{\it every\/} Latt\`es map with \[L'=a\] must have a postcritical fixed point.
However, when \[a\] is odd, the
period two orbits
$$ 0~\leftrightarrow~1/2~,\qquad \gamma/2~\leftrightarrow~(\gamma+1)/2 $$
in \[\T(\gamma)\] map to ramified period two orbits in \[\T(\gamma)/G_2\],
and there is no postcritical fixed point.

 Caution: In this last case, we can no longer
realize every conjugacy class of maps by restricting \[\gamma\]
to the Siegel region. A larger fundamental domain is needed.
For explicitly worked out examples in both cases,
see equations (15), (18) and (19) below; and
for further discussion see \S7.\ms

Here is another characterization.

{\QP{\bf Lemma 5.6.} \it A Latt\`es map is flexible if and only if
the multiplier for every periodic orbit is an integer.\ss}

{\bf Proof.} This follows from Corollary 3.9. If \[n>2\] or if
\[a\not\in\Z\], then we can find infinitely many integers \[p>0\] so that
\[\omega\,a^p\not\in\Z\] for some \[\omega\in G_n\]. The number of
fixed points of the map \[\omega\,L^{\circ p}\] on the torus \[\T\]
grows exponentially
with \[p\] (the precise number is \[|\omega\,a^p-1|^2\]), and
each of these maps to a periodic point of the associated Latt\`es map \[f\].
If we exclude the three or four postcritical points, then the derivative
of \[f^{\circ p}\] at such a point will be \[\omega\,a^p\], so that the
multiplier of this periodic orbit cannot be an integer.\QED\ss

It seems very likely that power maps, Chebyshev maps,
and flexible Latt\`es maps
are the only rational maps such that the multiplier of every periodic
orbit is an integer. (For a related result, see Lemma 7.1 below.)
\bs

{\bf\S6. Latt\`es Maps before Latt\`es.} Although the name
of Latt\`es has become
firmly attached to the construction studied in this paper, it actually occurs
much earlier in the mathematical literature. Ernst Schr\"oder, in a well
known 1871 paper, first described ``Chebyshev'' type examples using
trigonometric functions, and then
gave an explicit one-parameter family of ``Latt\`es'' type examples as follows.
Let  \[x=\sn(u)\] be the Jacobi sine function with elliptic modulus $k$,
defined by the equation
$$ u\=\int_0^{\,x} {d\xi\over{\sqrt{(1-\xi^2)(1-k^2\xi^2)}}}~. $$
More explicitly, for any parameter \[k^2\ne 0, 1\], let
\[E_k\subset \C^2\] be the elliptic curve defined
by the equation \[y^2=(1-x^2)(1-k^2x^2)\].
Then the holomorphic differential
\[dx/y\] is smooth and non-zero everywhere on \[E_k\] (even at the
two points at infinity in terms of suitable local coordinates).
The integrals \[\oint dx/y\] around
cycles in \[E_k\] generate a lattice \[\Lambda\subset\C\], and the integral
$$ u(x,y)=\int_{(0,1)}^{(x,y)} dx/y $$
along any path from \[(0,1)\] to \[(x,y)\] in \[E_k\] is well defined
modulo this lattice. In fact the resulting coordinate \[u\] parametrizes
the torus \[\T=\C/\Lambda\], and we can set \[x=\sn(u)\] and
\[y=\cn(u)\,\dn(u)\]. Here \[\sn(u)\] is the Jacobi sine function, and
\[\cn(u)\] and \[dn(u)\]
are closely related doubly periodic meromorphic functions which satisfy
$$  \cn^2(u)\=1-\sn^2(u)\qquad{\rm and}\qquad \dn^2(u)\=1-k^2\sn^2(u)~.$$
Furthermore
$$\sn(2u)\= {{2\,\sn(u)\,\cn(u)\,\dn(u)}\over{1-k^2\sn^4(u)}}~. $$
(Compare [WW, \S22.2].)
Setting \[z=x^2=\sn^2(u)\]
it follows easily that there is a well defined rational function
$$ f(z)\= {4z(1-z)(1-k^2z)\over (1-k^2z^2)^2} \eqno(15)$$
of degree four which satisfies the semiconjugacy relation
$$   \sn^2(2u)\=f\Big(\sn^2(u)\Big)~. $$
This is Schr\"oder's example (modulo a minor misprint). In the terminology of
\S5, \[f\] is a ``flexible Latt\`es map'', described 47 years before Latt\`es.

It is not hard to see that
\[\sn(u)\] has critical values \[\pm 1\] and \[\pm 1/k\], and hence that
\[\sn^2(u)\] has critical values \[1\,,~1/k^2\,,~0\], and \[\infty\].
On the other hand the map \[f\] has three critical values  \[1\,,~1/k^2\]
and \[\infty\], which all map to the fixed point \[0=f(0)\].
Each of these three critical values is the
image under \[f\] of two distinct simple critical points.\ss

Lucyan B\"ottcher in 1904 cited the same example (with a different version of
the misprint). He was perhaps the first to think of this example from a
dynamical viewpoint, 
and to use the term ``chaotic'' to describe the behavior of the sequence
of iterates of \[f\]. In fact he
described an orbit \[z_0\mapsto z_1\mapsto\cdots\] as {\bit
chaotic\/} if for every convergent subsequence \[\{z_{n_i}\}\] the differences
\[n_{i+1}-n_i\] are unbounded. (Note that this includes examples such as
irrational rotations which are not chaotic in the modern sense.)

B\"ottcher actually cited a much earlier paper, written by Charles Babbage
in 1815, for a fundamental property of what we now call semiconjugacies.
For example, in order
to find a periodic point \[\psi^nx=x\] of a mapping \[\psi\], Babbage proceeded
as follows (see [Ba, p.~412]):

{\QP``{\it Assume as before \[~\psi x=\phi^{-1}f\phi\, x\], then
$$ \eqalign{\psi^2x=\phi^{-1}f\phi\,\phi^{-1}f\phi\, x=\phi^{-1}f^2\phi\, x~~ \cr
\psi^3 x=\phi^{-1} f^2\phi\,\phi^{-1}f\phi\, x=\phi^{-1}f^3\phi\, x~, }$$
~~and generally \[~\psi^n x=\phi^{-1} f^n\phi\, x\], hence our equation
becomes\hfil\break \vskip -.1in
\[\qquad\qquad\qquad\qquad  \phi^{-1} f^n \phi\,x=x\].\qquad\[\cdots~~\]}''
\ms}

\ni In modern terminology, we would say that \[\phi\] is a
{\bit semiconjugacy\/} from \[\psi\] to \[f\].
It follows that any periodic point of \[\psi\] maps to a
periodic point of \[f\]; and furthermore (assuming that
\[\phi\] is finite-to-one and onto) any periodic
point of \[f\] is the image of a periodic point of \[\psi\]. B\"ottcher
pointed out that the use of such an intermediary map \[\phi\] to relate the
dynamic properties of two maps \[\psi\] and \[f\] lies at the heart of
Schr\"oder's example.

J. F. Ritt carried out many further developments of these ideas in the 1920's.
(For further historical information, see [A].)
\bs

{\bf\S7. More Recent Developments.} This concluding section will outline
some of the special properties shared by some or all finite quotients of affine
maps. 

We first consider the class of {\bit flexible\/} Latt\`es maps,
as described in 5.5 and 5.6. These
are the only known rational maps without attracting cycles
which admit a continuous family of deformations preserving the
topological conjugacy class. In fact the \[C^\infty\] conjugacy class
remains almost unchanged as we deform the torus. Differentiability
fails only at the postcritical points; and the multipliers of periodic
orbits remain unchanged even at these postcritical points.

Closely related is the following:

{\QP{\bf Fundamental Conjecture.} The flexible Latt\`es maps
are the only rational maps which admit an ``invariant line field'' on their
Julia set.\ss}

\ni By definition \[f\] has an {\bit invariant line field\/} if its
Julia set \[J\] has positive Lebesgue measure, and if there is a measurable
\[f$-invariant field of real one-dimensional subspaces of the tangent bundle
of \[\widehat\C\] restricted to \[J\].
The importance of this conjecture is demonstrated by the following.
(See [MSS], and compare the discussion in [Mc2] as well as [BM].)

{\QP{\bf Theorem of Ma\~n\'e, Sad and Sullivan.}  \it If this
Fundamental Conjecture is true, then hyperbolicity is dense among
rational maps. That is, every rational map can be
approximated by a hyperbolic map.\ss}

\ni To see that every flexible Latt\`es map has such an invariant line field,
note that any torus \[\C/(\Z\oplus\gamma\Z)\] is foliated by a
family of circles \[\Im(t)={\rm constant}\] which is invariant under the affine
map \[L\]. If \[f\]
is the associated Latt\`es map \[L/G_2\], then this circle foliation maps
 to an \[f$-invariant
foliation of \[J(f)=\widehat\C\] which is not only measurable but actually
smooth, except for singularities at the four postcritical points.\ss

Let us define the {\bit multiplier
spectrum\/} of a degree \[d\] rational map \[f\] to be the function which
assigns to each \[p\ge 1\] the unordered list of multipliers at the \[d^p+1\]
(not necessarily distinct)
fixed points of the iterate \[f^{\circ p}\]. Call two maps {\bit isospectral\/}
if they have the same multiplier spectrum.

{\QP{\bf Theorem of McMullen.} \it The flexible Latt\`es maps are the only
rational maps which admit non-trivial isospectral deformations.
The conjugacy class of any rational map which is not flexible Latt\`es
is determined, up to finitely many choices, by its multiplier spectrum.\ss}

This is proved in  [Mc1, \S2]. McMullen points out that the Latt\`es maps
\[L/G_2\] associated with imaginary quadratic number fields provide a
rich source of isospectral examples which are not flexible.
First note the following.

{\QP{\bf Lemma 7.1.} \it Two Latt\`es maps \[L/G_2:\T/G_2\to\T/G_2\]
are isospectral if and only if they have the same derivative \[L'=a\],
up to sign, and the same numbers of periodic orbits of
various periods within the postcritical set \[P_f\].\/\ss}

{\bf Proof Outline.} Let \[\omega=\pm 1\]. The number of fixed points of
the map
 \[\omega L^{\circ p}\]  on the
torus  can be computed as  \[|\omega a^p - 1|^2\]. These fixed points occur
in pairs \[\{\pm t\}\], and each such pair
corresponds to a single fixed point of
the corresponding iterate \[f^{\circ p}\], where \[f\cong L/G_2\].
Whenever \[t\ne -t\] on the torus, the multiplier of
\[f^{\circ p}\] at this fixed point is also equal to
\[\omega a^p\]. For the exceptional points \[t=-t\] which are fixed under the
action of \[G_2\] and correspond to postcritical points of \[f\],
the multiplier of \[f^{\circ p}\] is equal to  \[a^{2p}\]. Thus, to
determine the multiplier spectrum completely, we need only know how many
points of various periods there are in the postcritical set.\QED\ss

{\bf Example 7.2.} Let \[\xi=i\sqrt k\] where \[k>0\] is a
square-free integer, and let \[a=m\xi+n\].
Then for each divisor \[d\] of \[m\] the lattice \[\Z[\xi d]
\subset\C\] is a \[\Z[a]$-module. Hence the linear map \[L(\tau)=a\tau\] acts
on the associated torus \[\T/\Z[\xi d]\]. If \[m\] is highly divisible, then
there are many possible choices for \[d\]. Suppose, to simplify the discussion,
that \[mk\] and \[n\] are both even, so that \[a^2\] is divisible by two
in \[\Z[a]\]. Then multiplication by \[a^2\] acts as the zero map on the group
consisting of elements of order two in \[\T\]. Thus
\[0=L(0)\] is the only periodic point under this action, hence the image
of \[0\] in \[\T/G_2\] is the only postcritical periodic point of \[L/G_2\].
It then follows from 7.1 that these examples are all isospectral.\bs

Berteloot and Loeb [BL1] have characterized Latt\`es maps in terms of the
dynamics and geometry of the associated homogeneous polynomial map of
\[\C^2\]. Every rational map \[f:\P^1\to\P^1\] of degree two or more
lifts to a homogeneous
polynomial map \[F:\C^2\to\C^2\] of the same degree with the origin as an
attracting fixed point. They show
that \[f\] is Latt\`es if and only if the boundary of the basin of the origin
under \[F\] is smooth and strictly pseudoconvex within some open set.
In fact, this boundary must be spherical in suitable local holomorphic
coordinates, except over the postcritical set of \[f\].

Anna Zdunik [Z] has
characterized Latt\`es maps using only measure theoretic properties.
It is not hard to see that the
standard probability measure on the flat torus pushes forward under \[\Theta\]
to an ergodic
\[f$-invariant probability measure on the Riemann sphere. This measure
is smooth and in fact real analytic, except for mild singularities at the
postcritical points. Furthermore, it is {\bit balanced\/},
in the sense that the preimage \[f^{-1}(U)\] of any simply connected subset of
\[\widehat\C\ssm P_f\] is a union of \[n\] disjoint sets of equal measure.
Hence it
coincides with the unique measure of maximal entropy, as constructed by
Lyubich [Ly], or by Freire, Lopez and M\~an\'e [FLM], [Ma]. The converse
theorem is much more difficult:

{\QP{\bf Theorem of Zdunik.} \it The Latt\`es maps are the only rational maps
for which the measure of maximal entropy is absolutely continuous with respect
to Lebesgue measure.\ms}

We can think of the maximal entropy measure \[\mu_{\max}\] as describing
the asymptotic distribution of random backward orbits. That is, if we start
with any non-exceptional initial point \[z_0\], and then use a fair \[d$-sided
coin or die to iteratively choose a backward orbit
$$	\cdots \mapsto z_{-2}\mapsto z_{-1}\mapsto z_0~, $$
then \[\{z_n\}\] will be equidistributed with respect to \[\mu_{\max}\].
This measure \[\mu_{\max}\] always exists. An absolutely continuous
invariant measure is much harder to find, and an invariant measure which
is ergodic and
belongs to the same measure class as Lebesgue measure is even rarer.
However Latt\`es maps are not the
only examples: Mary Rees [Re] has shown that for
every degree \[d\ge 2\] the moduli space of degree \[d\] rational maps
has a subset of positive measure consisting of maps \[f\] which have an
ergodic invariant measure \[\mu\] in the same measure class as Lebesgue
measure. Such a measure is clearly unique, since Lebesgue almost every
forward orbit \[z_0\mapsto z_1\mapsto z_2\cdots\] must be equidistributed
with respect to \[\mu\].

Using these ideas, an easy consequence of Zdunik's Theorem is
the following.

 {\QP{\bf Corollary.} \it A  Lebesgue randomly chosen forward
orbit for a Latt\`es map has the same asymptotic distribution as a
randomly chosen backward orbit.\ss}

\ni I don't know whether Latt\`es maps are the only ones with this property.\ss

In general, different rational
maps have different invariant measures, except that every invariant measure for
\[f\] is also an invariant measure for its iterates \[f^{\circ p}\]. However,
every Latt\`es map \[L/G_n\] shares its measure \[\mu_{\max}\] with a rich
collection of Latt\`es maps \[\widehat L/G_n\] where \[\widehat L\] ranges
over all affine maps of the torus  which commute with the action of \[G_n\].
This collection forms a semigroup which is not finitely generated. (If we
consider only the linear torus maps \[\widehat L(\tau)=a\tau\], then we obtain
a commutative semigroup.) I don't know any other examples, outside of the
Chebyshev and power maps, of a semigroup of rational
maps which is not finitely generated, and which shares
a common non-atomic invariant measure. (See [LP] for related results.)\ss

Closely related is the study of commuting rational maps. Following a
terminology introduced much later
by Veselov [V], let us call a rational map \[f\]
{\bit integrable\/} if it commutes with another rational map, \[f\circ g=g\circ
f\], where both \[f\] and \[g\] have degree at least two, and where no iterate
of \[f\] is equal to an iterate of \[g\].

{\QP{\bf Theorem of Ritt and Eremenko.} \it A rational map \[f\] of degree
\[d_f\ge 2\] is integrable if and only if it is a finite quotient
of an affine map; that is if and only if it
is either a Latt\`es, Chebyshev, or power map. Furthermore, the commuting
map \[g\] must have the same Julia set, the same flat orbifold metric,
the same measure of maximal entropy, and the same set of preperiodic points
as \[f\].\ms}

This is a modern formulation of a statement which was proved by Ritt [R2]
in 1923, and by Eremenko [E] using a quite different method in 1989. For higher
dimensional analogues, see [DS]. In fact there has been a great deal of
interest in higher dimensional analogues of Latt\`es maps in recent years.
Compare [BL2], [Di], [Du], [V].\bs

{\bf\S8. Examples.} This concluding section will provide explicit formulas
for some particular Latt\`es maps.

{\bf 8.1. Degree Two Latt\`es Maps.} Recall from Lemma 5.4 and equation (14)
that the derivative \[L'=a\] for a torus map of degree \[d\] must
either be a (rational) integer, so that \[d=a^2\], or must be an imaginary
quadratic integer of the form
\[a=\Big(q\pm\sqrt{q^2-4d}\,\Big)/2\] with \[q^2<4d\], satisfying
\[a^2-q\,a+d=0\] and \[|a|^2=d\]. Furthermore, replacing
\[a\] by \[-a\] if necessary, we may assume that \[q\ge 0\]. Thus,
in the degree two case, the only distinct possibilities are \[q=0,1,2\], with
$$a\= i\sqrt 2~,\qquad{\rm or}\qquad a\=(1\pm i\sqrt 7)/2~,\qquad{\rm or}
\qquad a\=1\pm i~. $$
In each of these cases, the associated torus \[\T=\C/\Lambda\] is necessarily
conformally isomorphic to the quotient \[\C/\Z[a]\]. In fact we can
assume that \[\Lambda=\Z\oplus\gamma\Z\] with \[\gamma\] in the Siegel region
(8), and set \[a=r+s\gamma\] with \[r\,,\,s\in\Z\]. Let us assume, to fix
our ideas, that \[\Im(a)>0\]. Then, arguing as in the proof of 5.4, we have
$$  0~<~s~\le~{2\,\Im(a)\over\sqrt 3}~\le~{2|a|\over\sqrt 3}\=
{2\sqrt 2\over\sqrt 3}~\approx~ 1.63~.$$
Therefore \[s=1\], hence \[a\equiv\gamma~(\mod\Z)\]; so the lattice \[\Lambda\]
must be equal to \[\Z[a]\].

First suppose that \[f\cong L/G_2\] is a Latt\`es map of type \[\{2,2,2,2\}\],
with \[L(t)=at+b\]. The four
points of the form \[\lambda/2\] in \[\T\] map to the four postcritical
points of \[f\]. Hence the
action of the Latt\`es map \[f\] on its postcritical set is mimicked by
the action of \[L\] on this group of elements of the form \[\lambda/2\] in
\[\T\], or equivalently by the action of \[t\mapsto at+2b\] on
the four element group \[\Z[a]/2\Z[a]\].
 A brief computation shows that the quotient group
\[\Z[a]/\Big(2\Z[a]+(a-1)\Z[a]\Big)\] of Theorem 5.1 is trivial when \[q\] is
even but has two elements when \[q\] is odd. Thus, in the two cases
\[a=i\sqrt 2\] and \[a=1+i\] where \[q\] is even,
we may assume that \[L(t)=at\]. In these cases, the equation
\[a^2-qa+2=0\] implies that \[a^2\equiv 0~~\mod 2\Z[a]\], and hence that
$$ 1\mapsto a\mapsto 0\qquad{\rm and}\qquad 1+a\mapsto a\mapsto 0 $$
in this quotient group. Thus in these two cases there is a unique postcritical
fixed point, which must have mutiplier \[a^2\] by 3.9. In fact,
the diagram of critical and postcritical
points for the Latt\`es map \[f\] necessarily has the following form.
$$      \matrix{
          &       &         &             &\obul \cr
          &       &         &             &\mapsup\cr
* &\mapsto&\bullet&\mapsto&\bullet&\mapsfrom&\bullet&\mapsfrom & *}$$
Here each star stands for a simple critical point,
each heavy dot stands for a ramified (or postcritical) point, and the
heavy dot with a circle around it stands for a postcritical fixed point.
If we put the two critical points at \[\pm 1\] and put
the postcritical fixed point at infinity, then \[f\] will have the form
$$           f(z)\=(z+z^{-1})/a^2\,+\,c $$
for some constant \[c\]. In fact we easily derive the forms
$$\eqalign{  &f(z)\= -(z+z^{-1})/2+\sqrt 2\qquad{\rm when}\qquad a=i\sqrt 2~,
\quad{\rm and} \cr
  &f(z)\= \pm(z+z^{-1})/2i\qquad\qquad{\rm when}\qquad a=1\pm i~.}$$
On the other hand, for \[a=(1\pm i\sqrt 7)/2\], a similar argument shows
that there are two possible critical orbit diagrams, as follows. Either:
$$ \matrix{
   {*}&\mapsto&\bullet&\mapsto&\obul&&
  &  {\bf*}&\mapsto&\bullet&\mapsto&\obul}
$$
with two postcritical fixed points, or
$$  \matrix{
    *&\mapsto&\bullet&\mapsto&\bullet&\longleftrightarrow&\bullet&
       \mapsfrom&\bullet&\mapsfrom&*}
$$
with no postcritical fixed points. In the first case, if we put the
postcritical fixed points at zero and infinity, and another fixed point
at \[+1\], then the map takes the form
$$ f(z)\=z{z+a^2\over a^2z+1}~. $$
This commutes with the involution \[z\mapsto 1/z\], and we can take the
composition
$$z~\mapsto f(1/z)\=1/f(z)\= {az^2+1\over z(z+a^2)} $$
as the other Latt\`es map with the same value of
\[a\], but with \[\{0,\infty\}\] as postcritical period two orbit.
(See [M5, \S B.3] for further information on these maps.\footnote{$^6$}{\tp
Caution: In both [M3, 2000] and [M5], the term ``Latt\`es map'' was used with
a more
restricted meaning, allowing only maps of type \[\{2,2,2,2\}\] with \[n=2\].})

We can also ask for Latt\`es maps of degree two of the form \[L/G_n\] with
\[n>2\]. However, only the type \[\{2,4,4\}\] with
n=4 can occur, since, of the lattices
\[\Z[a]\] described above, only \[Z[1+i]=\Z[i]\] admits a rotation of order
greater than \[2\]. In fact the rotation \[t\mapsto it\] of the torus
\[\C/\Z[i]\] corresponds to an
involution \[z\mapsto -z\] which commutes with the associated
Latt\`es map \[f(z)=(z+z^{-1})/2i\].
To identify \[z\] with \[-z\], we can introduce the new variable \[w=z^2\]
and set
$$g(w)\=g(z^2)\=f(z)^2\=-(z^2+2+z^{-2})/4\= -(w+2+w^{-1})/4~. $$
Up to holomorphic conjugacy, this is the unique degree two Latt\`es map
of the form \[L/G_4\]. Its critical points \[\pm 1\] have orbit
\[ 1~\mapsto~ -1~\mapsto 0~\mapsto~\infty~~\lloop\], so that \[-1\] is both a
critical point and a critical value, yielding the following schematic diagram.
$$\matrix { * & \mapsto & * & \mapsto & \bullet & \mapsto & \obul\cr
&& 2 && 4 && 4 }$$
Here the ramification index is indicated underneath each ramified point.
Thus the map has type \[\{2,4,4\}\], as expected.
(Alternatively, the map \[z\mapsto 1-2/z^2\], with critical points zero and
infinity and with critical orbit
\[~0~\mapsto~\infty~\mapsto~ 1~\mapsto~ -1~\,\lloop\], could also be used as
a normal form for this same conjugacy class.)

\ss

{\bf 8.2. Degree Three.} If the torus map \[L(t)=at+b\] has degree \[|a|^2=3\],
then according to equation (14) we can write
\[a=\Big(q\pm\sqrt{q^2-12}\,\Big))/2\] for some
integer \[q\] with \[q^2<12\] or in other words \[|q|\le 3\]. I will try to
analyze only a single case, choosing \[q=0\] with \[a=i\sqrt 3\] so that
\[a^2=-3\].

For this choice \[a=i\sqrt 3\], setting \[a=r+s\gamma\] as in the proof of 5.4,
we find that \[s\] can be either one or two, and it follows easily
that there are exactly two essentially distinct tori
which admit an affine map \[L\] with derivative \[L'=a\].
We can choose either the hexagonally symmetric torus
\[\T=\C/\Z[\omega_6]=\C/\Z[(a+1)/2]\], or
its 2-fold covering torus \[\T'=\C/\Z[a]\].

For the torus\[\C/\Z[a]\], since there is no
\[G_3\] or \[G_4\] symmetry, we are necessarily in the
case \[n=2\]. A brief computation shows that the quotient group
\[\Z[a]/(2\Z[a]+(a-1)\Z[a])\] of Theorem
5.1 has two elements, so there are two
possible Latt\`es maps, corresponding to the two affine maps \[L(t)=at\] and
\[L(t)=at+1/2\]. The corresponding critical orbit diagrams have the form
$$\matrix{ *&\mapsto&\obul&\quad&*&\mapsto&\obul&\quad& *
\mapsto&\bullet&\leftrightarrow&\bullet&\mapsfrom& *}$$
with two postcritical fixed points, and
$$\matrix{ *&\mapsto&\bullet &\mapsto&\bullet&\mapsfrom& *\cr
&&\mapsup&&\mapsdown\cr
 *&\mapsto&\bullet &\mapsfrom&\bullet&\mapsfrom& *}$$
with no postcritical fixed points. In the first case, if we place the
postcritical fixed points at
zero and infinity, and place a fixed point with multiplier \[+a\] at \[+1\],
then the map takes the form
$$  f(z)\={z(z-a)^2\over(az-1)^2}~. $$
The remaining fixed point then lies at \[-1\] and has multiplier \[-a\].
The two remaining critical points \[\pm 2i-a\] map to the period two
postcritical orbit \[ -2i+a~\longleftrightarrow~2i+a~. \]

We can construct the other Latt\`es maps with the same \[a\] and the same
lattice \[\Z[a]\] by composing \[f\] with the M\"obius
involution \[g\] which satisfies
$$g~:~0~\leftrightarrow ~ 2i+a~,\qquad g~:~\infty~\leftrightarrow  ~-2i+a~. $$
The critical orbit diagram for this composition  permutes
the four postcritical points cyclically, as required.

{\bf A beautifully symmetric example.}
Now consider the torus \[\T=\C/\Z[\omega_6]\]. As noted at the end of \S4,
the quotient \[\T/G_2\], with its flat orbifold metric, is isometric
to a regular tetrahedron with the four cone points as vertices.
Again there are two distinct Latt\`es maps with invariant \[a^2=-3\]. The map
\[L(t)=at\] induces a highly symmetric piecewise linear map \[L/G_2\] of
this tetrahedron. (Compare [DMc] for a discussion of symmetric rational maps.)
The four vertices are postcritical fixed points of this map,
and the midpoints of the four faces are the critical points,
each mapping to the opposite vertex. Thus the critical
orbit diagram has the following form.
$$\matrix{ * & \mapsto & \obul & \quad & * & \mapsto & \obul & \quad &
 * & \mapsto & \obul & \quad & * & \mapsto & \obul }$$
Similarly,
the midpoint of each edge maps to the midpoint of the opposite edge, forming
a period two orbit.

If we place these critical points on the Riemann sphere at the
cube roots of \[-1\] and at infinity, then this map takes the form
$$	f(z)\={6\,z\over z^3-2}~, \eqno(16) $$
with a critical orbit \[\omega~\mapsto~-2\omega~\lloop\]
whenever \[\omega^3=-1\], and also
\[\infty~\mapsto~0~\lloop\].

The affine map \[L(t)=at+1/2\] yields a Latt\`es map \[L/G_2\]
with the same critical and
postcritical points, but with the following critical orbit diagram.
$$\matrix{ *&\mapsto&\bullet&\leftrightarrow&\bullet&\mapsfrom&*&\quad&
 *&\mapsto&\bullet&\leftrightarrow&\bullet&\mapsfrom&*} $$
Such a map can be constructed by composing the map \[f\] of (16) with
the M\"obius involution $$g(z)=(2-z)/( 1+z)$$ which satisfies
\[ ~ -1~\leftrightarrow~\infty~~{\rm and}~~ \omega_6~\leftrightarrow
~\overline\omega_6~.\] This corresponds to a $180^\circ$ rotation of the
tetrahedron about an axis joining the midpoints of two opposite faces.\ss

Now consider Latt\`es maps \[L/G_n\] with \[n\ge 3\] and with \[a=i\sqrt 3\].
Evidently the lattice must be \[\Z[\omega_6]\], and \[n\] must be either
3 or 6, so the type must be either \[\{3,3,3\}\] or \[\{2,3,6\}\].
Using 5.2, we can easily check that
there is just one possible map in each case, corresponding to the linear map
\[L(t)=at\]. Since both \[G_2\]
and \[G_3\] are subgroups of \[G_6\], this torus map \[L(t)=at\]
gives rise to maps of type \[\{2,2,2,2\}\] and \[\{3,3,3\}\] and
\[\{2,3,6\}\] which are related by the commutative diagram
$$	\matrix { L & \mapsto & L/G_2~~\cr
\mapsdown&&\mapsdown\cr L/G_3 & \mapsto& L/G_6~.} $$
Here \[L/G_2\] is the ``beautifully symmetric example'' of equation (16).
The corresponding Latt\`es map \[L/G_6\] of type \[\{2,3,6\}\]
can be constructed from (16) by identifying
each \[z\] with \[\omega\,z\] for \[\omega\in G_3\].
If we introduce the new variable \[\zeta=z^3\], then the corresponding
map \[L/G_6=f/G_3\] is given by mapping \[\zeta=z^3\] to \[g(\zeta)=f(z)^3\],
so that
$$ g(\zeta)\=
\left({6\,z\over z^3-2}\right)^3\={6^3\,\zeta\over(\zeta-2)^3}~. \eqno(17)$$
The three critical points at the cube roots of \[-1\] now coalesce into
a single critical point at \[-1\], with \[g(-1)=g(8)=8\].
There is still a critical point at infinity with \[g(\infty)=g(0)=0\].
But now infinity is also a critical value. In fact there if a double critical
point at \[\zeta=2\], with \[g(2)=\infty\]. The corresponding diagram for
the critical and postcritical points \[2\mapsto\infty\mapsto 0\] and
\[-1\mapsto 8\]  takes the form
$$\matrix{	** & \mapsto&  * & \mapsto&  \obul&\qquad\qquad &   *~  \mapsto&  \obul\cr
&&                3 &&     6 &&             & 2 } $$
where the symbol \[*\!*\] stands for a critical point of multiplicity two.
The multipliers at the two postcritical fixed points are \[a^6=-27\] and
\[a^2=-3\] respectively.\ss

Similarly we can study the Latt\`es map \[L/G_3\].
In this case the three points of \[\T\] which are fixed by \[G_3\]
all map to zero. Thus the three cone points of the orbifold
\[\T/G_3\] all map to one of the three. The corresponding diagram
has the following form.
$$\matrix{		**&  \mapsto&  \bullet&  \mapsto &  \obul&  \mapsfrom&
  \bullet&  \mapsfrom&  **\cr
 &&                    3 &&      3 &&     3 } $$
If we put the critical points at zero and infinity, and the postcritical fixed
point at \[+1\] (compare [M4]), then this map takes the form
$$	f(z) \= {z^3 + \omega_3\over \omega_3 z^3 + 1} ~, $$
with critical orbits   \[0 \mapsto \omega_3 \mapsto 1~\lloop\], and
\[\infty \mapsto \overline\omega_3 \mapsto 1~\lloop\].
In contrast to \[L/G_2\] and \[L/G_6\], this cannot be represented
as a map with real
coefficients. In fact the invariant \[a^3=-i\sqrt{27}\] is not a real number,
so this \[f\] is not holomorphically conjugate to
its complex conjugate or mirror image map.
(For a similar example with \[a^3=-8\] which occurs in the study of rational
maps of the projective plane, see [BDM, \S4 or \S6].)

Note that \[f\] commutes with the involution \[z\mapsto 1/z\]. If we
identify \[z\] with \[1/z\] by introducing a new variable \[w=z+1/z\], then
we obtain a different model for \[L/G_6\], which is necessarily
conformally conjugate to (17).\ms

{\bf 8.3. Flexible Latt\`es maps.} Recall from \S5 that there is just one
connected family of flexible Latt\`es maps of degree \[a^2\] for each even
integer \[a\], but that there are two distinct families of degree \[a^2\]
when \[a\] is odd. For \[a^2=4\],
the Schr\"oder family (15), constructed by expressing \[\sn^2(2t)\]
as a rational function of \[\sn^2(t)\], exhausts all of
the possibilities. This family depends on a parameter
\[k^2\in\C\ssm\{0,1\}\] and has postcritical set
\[\{0\,,\,1\,,\,\infty\,,\,1/k^2\}\], with all postcritical points mapping
to the fixed point zero. Using the corresponding formula for \[\sn^2(3t)\]
and following Schr\"oder's method, we obtain the family
$$ f(z)\={z\,(k^4z^4-6k^2z^2+4(k^2+1)z-3)^2\over
(3k^4z^4-4k^2(k^2+1)z^3+6k^2z^2-1)^2} \eqno (18) $$
of degree nine Latt\`es maps, with the same postcritical set
\[\{0\,,\,1\,,\,\infty\,,\,1/k^2\}\], but with all postcritical points
fixed by \[f\]. Note that \[f\] commutes with the involution
\[z\mapsto 1/k^2z\] which permutes the postcritical points. Composing \[f\]
with this involution, we obtain a different family
$$ z~\mapsto~ {1\over k^2f(z)}\=f\left({1\over k^2 z}\right)\=
{(3k^4z^4-4k^2(k^2+1)z^3+6k^2z^2-1)^2\over
 k^2 z((k^4z^4-6k^2z^2+4(k^2+1)z-3)^2} \eqno(19) $$
with the same postcritical set, but with all postcritical orbits of period two.
Higher degree examples could be worked out by the same method. Presumably
they look much like the degree four case for even degrees, and much
like the degree nine case for odd degrees.\ms

{\bf References.}\ss

\ref [A] D.~S.~Alexander, ``A History of Complex Dynamics'', Vieweg 1994.\ss

\ref [Ba] C.~Babbage. {\it An essay on the calculus of functions\/}, Phil.
Trans. Royal Soc. London {\bf 105} (1815) 389-423. (Available through
www.jstor.org .)\ss

\ref [BL1] F.~Berteloot and J.-J.~Loeb, {\it Spherical hypersurfaces and
Latt\`es rational maps\/},  J.~Math. Pures Appl. {\bf 77} (1998) 655-666.

\ref [BL2] F.~Berteloot and J.-J.~Loeb, {\it
Une caract\'erisation g\'eom\'etrique des exemples de Latt\`es de
${\Bbb P}\sp k$\/}, Bull. Soc. Math. France {\bf 129} (2001), 175--188.

\ref [BM] F.~Berteloot and V.~Mayer, ``Rudiments de dynamique holomorphe'',
Soc. Math. France 2001.

\ref [BDM] A.~Bonifant, M.~Dabija and J.~Milnor, {\it Elliptic curves as
attractors in \[{\bbb P}^2\], Part 1\/}, in preparation.\ss

\ref [B\"o]  L. E. B\"ottcher,
{\it The principal laws of convergence of iterates and
their application to analysis\/} (Russian), Izv. Kazan. Fiz.-Mat. Obshch.
{\bf 14} (1904) 155-234.

\ref [De] L. DeMarco, {\it Iteration at the boundary of the space of rational
maps\/}, Preprint, U. Chicago 2004.

\ref [Di] T.-C.~ Dinh, {\it Sur les applications de Latt\`es de
\[{\bbb P}^k\]\/}, J.~Math.~Pures Appl. {\bf 80} (2001) 577-592.

\ref [DS] T.-C.~ Dinh and N.~Sibony, {\it Sur les endomorphismes holomorphes
permutables\/}, Math. Annalen {\bf 324} (2002) 33-70.

\ref [DH] A.~Douady and J.~Hubbard, {\it A proof of Thurston's topological
characterization of rational functions\/},
Acta.~Math.~{\bf 171} (1993) 263-297.\ss

\ref [DMc] P.~Doyle and C.~McMullen, {\it Solving the quintic by iteration\/},
 Acta Math.~{\bf 163} (1989) 151--180.\ss

\ref [Du] C.~Dupont, {\it
Exemples de Latt\`es et domaines faiblement sph\'eriques de \[\C^n\]\/},
Manuscripta Math. {\bf 111} (2003) 357--378.\ss

\ref [E] A.~Eremenko, {\it
Some functional equations connected with iteration of rational functions\/},
Algebra i Analiz {\bf 1} (1989) 102-116 (Russian);
Leningrad Math. J. = St. Peterburg Math. J. {\bf 1} (1990) 905-919 (English).

\ref [FLM] A.~Freire, A.~Lopes, and R.~Ma\~n\'e, {\it
An invariant measure for rational maps\/},
Bol. Soc. Brasil. Mat. {\bf 14} (1983), 45-62.\ss

\ref [K] G.~K{\oe}nigs, {\it Recherches sur les integrals de certains
equations fonctionelles\/}, Ann.~Sci.~\'Ec. Norm.~Sup. {\bf 1} (1884)
suppl.~1-41.\ss

\ref [La] S.~Latt\`es, {\it Sur l'iteration des substitutions rationelles
et les fonctions de Poincar\'e\/}, C.~R. Acad.~Sci.~Paris {\bf 166} (1918)
26-28. (See also ibid.~pp.~486-89.)\ss

\ref [LP] G.~Levin and F.~Przytycki, {\it When do two rational functions have
the same Julia set?\/}, Proc.~A.~M.~S. {\bf 125} (1997) 2179-2190.

\ref [Ly] M.~Lyubich, {\it Entropy properties of rational endomorphisms
of the Riemann sphere\/}, Ergodic Th.~and Dy.~Sys. {\bf 3} (1983) 351-385.\ss 

\ref [Ma] R.~Ma\~n\'e, {\it On the uniqueness of the maximizing measure for
rational maps.\/}  Bol. Soc. Brasil. Mat. {\bf 14}  (1983) 27--43.

\ref [MSS] R.~Ma\~n\'e, P.~Sad, and D.~Sullivan, {\it
On the dynamics of rational maps\/},
Ann.~Sci.~\'Ecole Norm.~Sup. {\bf 16} (1983) 193--217.

\ref [Mc1] C.~McMullen, {\it Families of rational maps and iterative
root-finding algorithms\/}, Ann. of Math. {\bf 125} (1987) 467--493.

\ref [Mc2] ---------, {\it Frontiers in complex dynamics\/}, Bull.~A.~M.~S.
{\bf 31} (1994) 155-172.

\ref [M1] J.~Milnor, {\it On the 3-dimensional Brieskorn manifolds
\[M(p,q,r)\]\/},
in ``Knots, Groups and 3-manifolds, Annals Math.~Studies {\bf 84}, Princeton
Univ.~Press 1975, or in Milnor  ``Collected Papers {\bf 2}'', Publish or
Perish 1995.\ss

\ref [M2] ---------, {\it Geometry and dynamics of quadratic rational maps\/},
       Experimental Math.~{\bf 2} (1993) 37-83.\ss

\ref [M3] ---------, ``Dynamics in One Complex Variable'',
 Vieweg, 1999, 2000; Princeton University Press 2004 (to appear).\ss

\ref [M4] ---------, {\it Rational maps with two critical points\/},
   Experimental Math.~{\bf 9} (2000) 481-522.\ss

\ref [M5] ---------, {\it Pasting together Julia sets, a worked out
example of mating\/}, Experimental Math.~{\bf 13} (2004) 55-92.\ss

\ref [P] H.~Poincar\'e, {\it Sur une classe nouvelle de transcendantes
uniformes\/}, J.~Math.~Pures Appl.~{\bf 6} (1890) 313-365 ({\OE}vres IV,
537-582).\ss

\ref [Re] M.~Rees, {\it
Positive measure sets of ergodic rational maps.\/},
Ann. Sci. \'Ecole Norm.~Sup. {\bf 19} (1986) 383--407.\ss

\ref [R1] J.~F.~Ritt, {\it Periodic functions with a multiplication theorem\/},
Trans. Amer. Math. Soc. {\bf 23} (1922) 16-25.

\ref [R2] ---------, {\it Permutable rational functions\/},
Trans. Amer. Math. Soc. {\bf 25} (1923) 399--448.\ss

\ref [R3] ---------, {\it Meromorphic functions with addition or
multiplication theorems\/}, Trans. Amer. Math. Soc. 29 (1927), 341-360.\ss

\ref [S] E.~Schr\"oder, {\it
Ueber iterirte Functionen\/}, Math. Ann. {\bf 3} (1871), 296-322. Available
through http://gdz.sub.uni-goettingen.de/en (search by title).\ss

\ref [V] A.~P.~Veselov, {\it Integrable maps\/}, Uspekhi Mat. Nauk {\bf 46:5}
(1991) 3-45; Russ. Math. Surveys  {\bf 46:5} (1991) 1-51.

\ref [WW] E.~T.~Whittaker and G.~N.~Watson, ``A Course of Modern Analysis'',
Cambridge U.~Press 1966 (4th edition).\ss

\ref [Z] A.~Zdunik, {\it Parabolic orbifolds and the dimension of the maximal
measure for rational maps\/}, Inv.~Math.~{\bf 99} (1990) 627-649.

\hskip 3.5in September 2004

\hskip 3.5in Institute for Mathematical Sciences

\hskip 3.5in Stony Brook University

\hskip 3.5in Stony Brook, NY 11794-3660 USA\medskip

\hskip 3.5in jack@math.sunysb.edu

\hskip 3.5in www.math.sunysb.edu/$\sim$jack

\end